\documentclass[12pt,a4paper]{amsart}

\usepackage[utf8]{inputenc}

\usepackage{mathrsfs}
\usepackage{todonotes}

\usepackage{amssymb, tikz}

\usepackage{hyperref}

\usepackage{amsmath}
\usepackage[utf8]{inputenc}

\newtheorem{theorem}{Theorem}[section]
\newtheorem{claim}[theorem]{Claim}

\newtheorem{lemma}[theorem]{Lemma}

\newtheorem{corollary}[theorem]{Corollary}

\theoremstyle{definition}
\newtheorem{definition}[theorem]{Definition}

\newtheorem{problem}[theorem]{Problem}

\theoremstyle{remark}
\newtheorem{remark}[theorem]{Remark}
\newtheorem{question}[theorem]{Question}
\newtheorem{notation}[theorem]{Notation}

\newcommand{\lex}{{\rm lex}}

\newcommand{\ZF}{{\rm ZF}}
\newcommand{\ZFC}{{\rm ZFC}}
\newcommand{\DC}{{\rm DC}}

\newcommand{\Borel}{{\rm Borel}}
\newcommand{\Ord}{{\rm Ord}}

\newcommand{\nor}{{\rm nor}}
\newcommand{\tr}{{\rm tr}}
\newcommand{\rk}{{\rm rk}}

\newcommand{\id}{{\rm id}}

\newcommand{\Suc}{{\rm Suc}}

\newcommand{\dom}{{\rm dom}}

\newcommand{\rng}{{\rm range}}

\newcount\skewfactor
\def\mathunderaccent#1#2 {\let\theaccent#1\skewfactor#2
\mathpalette\putaccentunder}
\def\putaccentunder#1#2{\oalign{$#1#2$\crcr\hidewidth
\vbox to.2ex{\hbox{$#1\skew\skewfactor\theaccent{}$}\vss}\hidewidth}}
\def\name{\mathunderaccent\tilde-3 }

\begin{document}

\title[ On the classification of definable ccc forcing notions]{On the classification of definable ccc forcing notions}

\author[M.  Golshani]{Mohammad Golshani}

\address{Mohammad Golshani, School of Mathematics, Institute for Research in Fundamental Sciences (IPM), P.O.\ Box:
	19395--5746, Tehran, Iran.}

\email{golshani.m@gmail.com}
\urladdr{http://math.ipm.ac.ir/~golshani/}

\author[H. Horowitz]{Haim Horowitz}

\address{Haim Horowitz, Einstein Institute of Mathematics
Edmond J. Safra Campus,
The Heebrew University of Jerusalem.
Givat Ram, Jerusalem, 91904, Israel.}

\email{haim.horowitz@mail.huji.ac.il}

\author[S. Shelah] {Saharon Shelah}
\address{Einstein Institute of Mathematics\\
Edmond J. Safra Campus, Givat Ram\\
The Hebrew University of Jerusalem\\
Jerusalem, 91904, Israel\\
 and \\
 Department of Mathematics\\
 Hill Center - Busch Campus \\
 Rutgers, The State University of New Jersey \\
 110 Frelinghuysen Road \\
 Piscataway, NJ 08854-8019 USA}
\email{shelah@math.huji.ac.il}
\urladdr{http://shelah.logic.at}
\thanks{ The first author's research has been supported by a grant from IPM (No. 1401030417). The
	third author's research has been partially supported by Israel Science Foundation (ISF) grant no:
	1838/19. This is publication 1097 of third author. }




\begin{abstract}
We show that for a Suslin ccc forcing notion $\mathbb Q$
adding a Hechler real, ``$\ZF+\DC_{\omega_1}+$all sets of reals are
$I_{\mathbb Q,\aleph_0}$-measurable'' implies the existence of
an inner model with a measurable cardinal. We also introduce a wide class of
Suslin ccc forcing notions which add a Hechler real, so that the above result applies to them.
\end{abstract}

\maketitle
\numberwithin{equation}{section}
\section{introduction}
This paper can be seen as part of a line of research motivated by the
following general problem:

\begin{problem}
 Classify the nicely definable forcing notions.
\end{problem}
Under $\ZFC$, definable is usually interpreted to
mean the class of Suslin posets, i.e. posets $\mathbb{P}$ such that the domain of $\mathbb{P}$ is an
analytic set of reals, and both the order and the incompatibility relation of $\mathbb{P}$
are analytic. See \cite{judah} for more discussion on Suslin forcing notions.

There are several ways to make the question more precise. An old question of Prikry asks if it is consistent that
every non-trivial ccc forcing notion adds a Cohen real or a random real?
By a theorem of Shelah \cite{Sh:480}, see also Velickovic \cite{vel}, a modified version of Prikry’s conjecture holds in the context of Suslin ccc forcing, i.e., every Suslin ccc forcing either adds a Cohen real or is a Maharam algebra. Indeed, Shelah  showed that any Suslin ccc forcing which is not $\omega^\omega$-bounding adds a Cohen real. Also in \cite{blasz}, Blaszyck and Shelah
showed that it is relatively consistent with ZFC that every nonatomic $\sigma$-centered
forcing notion adds a Cohen real.
For further discussion of the problem, see \cite{Sh:666}.

Given an infinite cardinal $\kappa$ and a tree-like partial order $\mathbb Q$, whose conditions are subtrees of $\omega^{<\omega},$  we can assign to the pair $(\mathbb Q, \kappa)$ an ideal $I_{\mathbb Q, \kappa}$
which is defined as the closure of
\[
\{X\subseteq \omega^{\omega} : (\forall p\in \mathbb Q)(\exists p\leq q)\big(lim(q) \cap X=\emptyset \big)\},
\]
under $\leq \kappa$ unions.
We also say a set of reals $X$ is
$I$-measurable, where $I$ is an ideal on $\omega^{\omega},$ if there exists a Borel set $B$ such that $X\Delta B \in I$, where
$\Delta$ denotes the symmetric difference operation.

By the celebrated result of Solovay \cite{solovay}, the existence of an inaccessible cardinal implies the consistency of
``$\ZF+\DC+$all sets of reals are $I_{\mathbb Q,\aleph_0}$-measurable'', where $\mathbb{Q}$
is one of the the random or Cohen forcing notions. By later results of Shelah \cite{Sh:176},
the existence of an inaccessible cardinal is needed to get the above result for random forcing, while it is not needed
for the case of Cohen forcing.
Also, by \cite{Sh:176}, $\DC_{\omega_1}$ implies the existence of a non-Lebesgue
measurable set, and hence ``$\ZF+\DC_{\omega_1}+$all sets of reals
are $I_{\mathbb Q,\aleph_0}$-measurable'' is inconsistent, when
$\mathbb Q$ is the random real forcing.
A longstanding open question of Woodin asks the following.
\begin{question}
\label{woodin} Is the theory ````$\ZF+\DC_{\omega_1}+$all sets of reals are $I_{\mathbb C,\aleph_0}$-measurable'' consistent, where $\mathbb C$
is the Cohen forcing.
\end{question}

Motivated by the  results of \cite{HwSh:1067} and \cite{HwSh:1094}, we discuss
the following variant of the above problems, which asks to classify
nicely definable ccc forcing notions based on the consistency strength of some regularity properties that one can derive from them.
\begin{problem}
Classify the Suslin ccc forcing notions according
to the consistency strength of  ``$\ZF+\DC_{\omega_1}+$all sets of reals are $I_{\mathbb Q,\kappa}$-measurable'',
where $\kappa \in \{\aleph_0,\aleph_1\}$.
\end{problem}
\begin{remark}
 We can replace the theory $\ZF+\DC_{\omega_1}$ with some other theories like $\ZF, \ZF+\text{AC}_{\omega}, \ZF+\DC, \ZF+\DC_{\omega_1}, \ZFC$ or other similar theories.
\end{remark}

In the first main result of the paper, we discuss the $I_{\mathbb Q,\aleph_0}$-measurability for  Suslin ccc forcing
notions which add a Hechler real, and prove the following.
\begin{theorem}
\label{t1}
Let $\mathbb Q$ be a Suslin ccc forcing
notion which adds a Hechler real.
Then the consistency of ``$\ZF+\DC_{\omega_1}+$every set of reals is $I_{\mathbb Q,\aleph_0}$-measurable''
implies the existence of an inner model of $\ZFC$ with a measurable cardinal.
\end{theorem}
Indeed the above theorem will be a special case of a more general theorem, see Section \ref{regularity}.
The above
result is extended  to other Suslin ccc forcing notions not adding Hechler reals
in a subsequent paper \cite{F1561} by the last two authors.

We then introduce a wide class of Suslin ccc forcing
notions which add a Hechler real. Let
$\mathbb{Q}_{\bold n}^1$ be the forcing notion defined in
 \cite{HwSh:1067} (see Section \ref{prel} for the definition). Then we prove the following.
\begin{theorem}
\label{t2}
Forcing with $\mathbb{Q}_{\bold n}^1$ adds a Hechler real.
\end{theorem}
In particular, it follows that the consistency of ``$\ZF+\DC_{\omega_1}+$every set of reals is $I_{\mathbb{Q}_{\bold n}^1,\aleph_0}$-measurable''
implies the existence of an inner model with a measurable cardinal.

We may also note that by \cite{HwSh:1067}, the theory ``$\ZF+$every set of reals is $I_{\mathbb{Q}_{\bold n}^1,\aleph_1}$-measurable$+$
there exists an $\omega_1-$sequence of distinct reals'' is consistent relative to $\ZFC$.
The problem
of finding forcing notions $\mathbb Q$ for which $\ZF+\DC_{\omega_1}+I_{\mathbb Q,\aleph_0}$-measurability
is consistent (maybe relative to large cardinals) remains open. In
 \cite{F1424},  the following result is proved, where
 $\mathbb{Q}_{\bold n}^2$ is the forcing notion introduced in
\cite{HwSh:1067}.
\begin{theorem} (\cite{F1424}) Suppose there is a measurable cardinal.
Then in a suitable generic extension, there is an inner model of ``$\ZF+\DC_{\omega_1}+$all
sets of reals are $I_{\mathbb{Q}_{\bold n}^2,\aleph_1}-$measurable''.
\end{theorem}
We may note that by \cite{HwSh:1067}, the forcing notion $\mathbb{Q}_{\bold n}^2$
adds no dominating reals (and hence no Hechler reals).

Our notation is standard. For a forcing notion $\mathbb{P}$ and two conditions $p, q$
in it, we use the notation $p \leq q$ to mean that $q$ is a stronger condition than $p$.

\section{The additivity of the ideals derived from a Suslin
ccc forcing notion adding a Hechler real}
\label{additivity}

In this section we show  that under $\ZF+\DC_{\omega_1}$, if $\mathbb Q$ is a Suslin ccc forcing notion
adding a Hechler real, then the additivity of $I_{\mathbb Q,\aleph_0}$
is $\aleph_1$. This will allow us to prove in the next section that
$\ZF+\DC_{\omega_1}+$measurability for the ideal derived from such
forcing notions, implies the existence of an inner model with a measurable
cardinal. A main concept in the following proof is a variant of the
rank function for Hechler forcing originally introduced in \cite{GiSh:412}.

Before we continue, let us recall the Helcher forcing.
\begin{definition}
\label{hechler} The Hechler forcing
 $\mathbb D$ is defined as follows:
 \begin{enumerate}
 \item A condition in $\mathbb{D}$ is a pair $p=(t_p, f_p),$ where
 \begin{enumerate}
 \item $t_p \in \omega^{<\omega},$
 \item $f_p \in \omega^\omega.$
 \end{enumerate}
\item Let $p, q \in \mathbb D.$ Then $p \leq q$  iff
 \begin{enumerate}
 \item $t_p \subseteq t_q,$
 \item $f_p \leq f_q$,
 \item $(\forall n \in \dom(t_q)\setminus \dom(t_p)) \big( t_q(n) \geq f_p(n) \big).$
 \end{enumerate}
 \end{enumerate}
\end{definition}
\begin{notation}
Given a condition $p=(t_p, f_p) \in \mathbb D,$ we call $t_p$ the trunk of $p$ and denote it by $\tr(p).$
\end{notation}
We may note that forcing with $\mathbb{D}$ adds a canonical name
 $$\name{\eta}_{dom}=\bigcup \{\tr(p): p \in \dot{G}_{\mathbb D}    \}\footnote{$\dot{G}_{\mathbb D}$ is the canonical $\mathbb D$-name for the generic filter.}$$
 for an element of $\omega^\omega$ which dominates every ground model function in
 $\omega^\omega.$
  The following definition is a variant of
 \cite[Definition 4.4.2.]{GiSh:412}.

\begin{definition}
\label{d-rank} Suppose $p \in \mathbb D$ and  $I=\{r_k : k<\omega\}$
is a  maximal antichain above $p$. Let also $A=\{\tr(r_k) : k<\omega\}$.
For every $t \in \omega^{<\omega}$ with
$\tr(p) \unlhd t$, we define $\rk_{p,A}(t) \in \Ord \cup \{ \infty\}$  by defining when $\alpha \leq \rk_{p,A}(t)$:
\begin{itemize}
\item $0 \leq \rk_{p,A}(t)$ is always true,

\item  $1 \leq \rk_{p,A}(t)$ iff for every $l\in [\ell g(\tr(p)),\ell g(t))$,
$f_p(l) \leq t(l)$, and there is no $s \in A$ such that:
\begin{itemize}
\item $p \leq (s,f_p \restriction [\ell g(s),\omega))$,
\item $\ell g(s) \leq \ell g(t)$, and
\item$l\in [\ell g(\tr(p)),\ell g(s)) \rightarrow t(l) \leq s(l)$.
 \end{itemize}
\item Suppose $\alpha>1$. Then  $\alpha \leq \rk_{p,A}(t)$ iff for every $\beta<\alpha$,
for infinitely many $k$, $\beta \leq \rk_{p,A}(t^{\frown}\langle k \rangle)$.
\end{itemize}
\end{definition}
The following easy lemma will be used later in this section.
\begin{lemma}
\label{rank}
Assume $p \in \mathbb D$ and $\tr(p)\unlhd t \in \omega^{<\omega}$. Let   $I=\{r_k : k<\omega\}$
be a  maximal antichain above $p$, and set $A=\{\tr(r_k) : k<\omega\}$.
Then
\begin{enumerate}
\item[(a)] If $\omega_1 \leq \rk_{p,A}(t)$, then $rk_{p,A}(t)=\infty$.

\item[(b)]  $\rk_{p,A}(\tr(p))<\omega_1$.
\end{enumerate}
\end{lemma}
\begin{proof}
(a). Suppose $\omega_1 \leq \rk_{p,A}(t)$. We  prove by induction on $\omega_1 \leq \epsilon$
that  $\epsilon \leq \rk_{p,A}(t)$.
For $\epsilon=\omega_1$ the claim follows from the assumption. Now suppose that $\epsilon>\omega_1$
and the lemma holds for all ordinals less than $\epsilon.$
For each $k<\omega$ set $\zeta_k=\rk_{p,A}(t^{\frown}\langle k \rangle)$.

If there exists $\zeta <\omega_1$ such that $\{k : \zeta<\zeta_k\}$
is finite, then
by Definition \ref{d-rank},
 $\rk_{p,A}(t) \leq \zeta+1<\omega_1$,
which is a contradiction.
Thus we may assume that for each $\zeta <\omega_1$, the set $A_\zeta=\{k : \zeta<\zeta_k\}$
is infinite.  Let
$$\zeta_*=\sup\{\zeta_k : k<\omega, \zeta_k <\omega_1\}.$$
Then $\zeta_* < \omega_1.$
By our assumption, $A_{\zeta_*}=\{k : \zeta_k >\zeta_*\}$ is infinite, therefore by the
choice of $\zeta_*$,
for each
$k\in A_{\zeta_*}$, we have $\omega_1 \leq \rk_{p,A}(t^{\frown}\langle k \rangle).$
Thus by the induction hypothesis
 $$\bigwedge_{\zeta<\epsilon} \bigwedge_{k\in A_{\zeta_*}} \zeta \leq \rk_{p,A}(t^{\frown}\langle k \rangle).$$
It now follows from Definition \ref{d-rank} that $\epsilon \leq \rk_{p,A}(t)$.

(b). Suppose on the contrary that  $\rk_{p,A}(\tr(p))\geq \omega_1$. Then by clause (a),  $\rk_{p,A}(\tr(p))=\infty,$ and hence, by
induction on $n<\omega$, we can
choose $t_n \in \omega^{\ell g(\tr(p))+n}$
such that:
\begin{enumerate}
\item $t_0=\tr(p)$,
\item  $\omega_1 \leq \rk_{p,A}(t_n)$,
\item  $m<n \rightarrow t_m \lhd t_n$.
\end{enumerate}
Let
$f=\bigcup_{n<\omega} t_n$ and
 $q=(\tr(p), f) \in \mathbb D$. We first show that $q$ extends the condition $p$. To see this, note that for every $n<\omega$,
 $1 < \rk_{p,A}(t_n)$,
and therefore by Definition \ref{d-rank}, $f_p(l) \leq t_n(l)$
for every $l\in [\ell g(\tr(p)), \ell g(t_n)]$, from which the result follows.

Since $I$ is a maximal antichain above $p$ and $p \leq  q$, we can find
 some $k<\omega$ such that $q$ and $r_k$ are compatible, in particular we have
 \begin{center}
$(*)$\quad\quad\quad  $l\in [\ell g(\tr(p)), \ell g(\tr(r_k))) \implies \tr(r_k)(l) \geq f(l)$.
 \end{center}
 Let $n> \ell g(\tr(r_k))$. Then for each $l$ as in $(*)$, $f(l)=t_n(l).$
By the choice of $A$,  $\tr(r_k) \in A$, thus the assumption $1\leq \rk_{p,A}(t_n)$
implies that
 \begin{center}
$(**)$\quad\quad $\exists l\in [\ell g(\tr(p)), \ell g(\tr(r_k))) \big( \tr(r_k)(l)<t_n(l)=f(l)\big)$.
 \end{center}
Putting $(*)$ and $(**)$ together, we get the desired contradiction.
\end{proof}

The main result in this section is the following.
\begin{lemma}
\label{positive} ($\ZF+\DC$)
\begin{enumerate}
\item[(a)] Let $M$ be an inner model of $\ZFC$ with
$\aleph_1^M=\aleph_1$,
and let $I=I_{\mathbb D, \aleph_0}$.
 Then the union of all Borel $I$-small sets coded in the model $M$ is $I$-positive.

 \item[(b)] Assume $\ZF+\DC_{\omega_1}$. Then $\text{add}(I_{\mathbb D, \aleph_0})$, the additivity of the ideal $I_{\mathbb D, \aleph_0}$, is equal to $\aleph_1.$
 \end{enumerate}
\end{lemma}

\begin{remark}
 The assumption of the the existence of an inner model $M$ of $\ZFC$ with $\aleph_1^M=\aleph_1$ follow from $\ZF+\DC_{\omega_1}$. To see this, note that the assumption
   $\ZF+\DC_{\omega_1}$ the existence of
a subset  $A\subseteq \omega_1$ such that $\aleph_1^{L[A]}=\aleph_1.$ Then
 $M=L[A]$ is an inner model of $\ZFC$ with
$\aleph_1^M=\aleph_1$.
\end{remark}

\begin{proof}
First note that clause (b) follows from clause (a), the above remark and absoluteness arguments.
In order to prove clause (a) of the lemma, we introduce a family $\langle  B_\epsilon: \epsilon < \omega_1 \rangle \in M$ of Borel sets
$B_\epsilon \in I$ such that their union $\bigcup_{\epsilon<\omega_1}B_\epsilon$ is not in $I$. To this end, we first define a sequence
of pairs $\langle (\Lambda_{\epsilon},h_{\epsilon}) : \epsilon<\omega_1 \rangle$
in $M$ as follows.

For $\epsilon<\omega_1$, let $Y_{\epsilon}$
be the set of pairs $(\Lambda,h)$
such that:
\begin{enumerate}
\item $\Lambda \subseteq \omega^{<\omega}$ is a tree such that:

\begin{enumerate}
\item if $t \in \Lambda$, then either  $\Suc_{\Lambda}(t)=\omega,$
or $\Suc_{\Lambda}(t)=\emptyset$,

\item if $t_1,t_2 \in \omega^k$, $t_1 \in \Lambda$
and $t_1(l) \leq t_2(l)$ for every $l<k$, then $t_2 \in \Lambda$,

\item $\Lambda$ has no infinite branches.
\end{enumerate}
\item $h: \Lambda \rightarrow \epsilon+1$ is a function such that:
\begin{enumerate}
\item $h(\langle \rangle)=\epsilon$,

\item if $t_1\lhd t_2$ are in $\Lambda$, then $h(t_1)> h(t_2)$,

\item if $h(t)=\zeta+1$ then for all $k<\omega,~ h(t^{\frown} \langle k \rangle)=\zeta$.

\item if $h(t)=\zeta$ where $\zeta$ is a limit ordinal,
then $\zeta \leq \lim_{k<\omega}h(t^{\frown} \langle k \rangle)$.
\end{enumerate}
\end{enumerate}

\begin{claim}
\begin{enumerate}
\item[(a)] $Y_{\epsilon} \neq \emptyset$ for every $\epsilon<\omega_1$.

\item[(b)] For every $\epsilon<\zeta<\omega_1$ and $(\Lambda_1,h_1) \in Y_{\epsilon}$,
there exists $(\Lambda,h) \in Y_{\zeta}$ such that $\Lambda_1 \subseteq \Lambda$.
\end{enumerate}
\end{claim}
\begin{proof}
{\bf Proof of (a).}
We prove the claim by induction on $\epsilon<\omega_1$.
This is clear for $\epsilon=0.$
Now suppose that $\epsilon=\zeta+1$ and the claim holds for $\zeta.$
Let $(\Lambda,h) \in Y_{\zeta}$ and define $(\Lambda', h')$
as follows:
\begin{itemize}
\item $\Lambda' = \bigcup_{k<\omega} \langle k \rangle ^{\frown} \Lambda,$
where
$\langle k \rangle ^{\frown} \Lambda = \{ \langle \rangle\} \cup \{ \langle k \rangle ^{\frown} t: t \in \Lambda \},$
\item $\dom(h')=\Lambda',$
\item $h'(\langle\rangle)=\epsilon$,
\item for $k<\omega$ and $t \in \Lambda$,
$h'(\langle k \rangle ^{\frown} t)=h(t).$
\end{itemize}
It is easily seen that $(\Lambda', h') \in Y_{\epsilon}.$
 Now suppose that
$\epsilon$ is a limit ordinal, and the claim holds for all $\zeta < \epsilon.$ Let $(\zeta_{k} : k<\omega)$
be an increasing sequence of ordinals cofinal in
 $\epsilon$. By the induction hypothesis,  choose an increasing
sequence$((\Lambda_{k},h_{k}) : k<\omega)$
such that $(\Lambda_{k},h_{k}) \in Y_{\zeta_{k}}$.
Define
 $(\Lambda', h')$
as follows:
\begin{itemize}
\item $\Lambda' = \bigcup_{k<\omega} \langle k \rangle ^{\frown} \Lambda_k,$

\item $\dom(h')=\Lambda',$
\item $h'(\langle\rangle)=\epsilon$,
\item for $k<\omega$ and $t \in \Lambda_k$,
$h'(\langle k \rangle ^{\frown} t)=h_k(t).$
\end{itemize}
Again, it is easy to show that $(\Lambda', h') \in Y_{\epsilon}.$

{\bf Proof of (b).}
Assume  $\epsilon<\zeta<\omega_1$ and $(\Lambda_1,h_1) \in Y_{\epsilon}$. By clause (a), pick some
 $(\Lambda_2,h_2) \in Y_{\zeta}$ and define $(\Lambda, h)=(\Lambda_1,h_1)+(\Lambda_2,h_2)$ as follows:
\begin{itemize}
\item $\Lambda=\Lambda_1 \cup \Lambda_2$.

\item $h: \Lambda  \rightarrow \zeta+1$ is defined as:
\begin{center}
 $h(\eta)=$ $\left\{
\begin{array}{l}
         h_1(\eta) \hspace{3.8cm} \text{ if } \eta \in \Lambda_1 \setminus \Lambda_{2};\\
         h_2(\eta) \hspace{3.8cm} \text{ if } \eta \in \Lambda_2 \setminus \Lambda_{1};\\
         \max\{h_1(\eta),h_2(\eta)\}  \hspace{1.5cm} \text{ if } \eta \in \Lambda_1 \cap \Lambda_2.
     \end{array} \right.$
\end{center}
\end{itemize}
It is easy to see that $(\Lambda,h) \in Y_{\zeta}$ and $\Lambda_1 \subseteq \Lambda$.
\end{proof}
Now fix a sequence $\langle (\Lambda_{\epsilon},h_{\epsilon}) : \epsilon<\omega_1 \rangle \in M$
such that $(\Lambda_{\epsilon},h_{\epsilon}) \in Y_{\epsilon}$.

Given a subtree $\Lambda \subseteq \omega^{<\omega}$ set
\[
\max(\Lambda)=\{t \in \Lambda: \Suc_\Lambda(t)=\emptyset   \}.
\]
For $\epsilon<\omega_1$,  $k<\omega$ and  $\Lambda =\Lambda_\epsilon$, we  define the following objects:

\begin{enumerate}
\item[$(*)_1$] $\Omega_{\Lambda, k}=\{t_0 ^{\frown} \langle 2n_0+1\rangle ^{\frown} t_1 ^{\frown} \langle 2n_1+1 \rangle ^{\frown} \cdots ^{\frown} t_k : n_i <\omega$
and $t_i \in \max(\Lambda) \}$.

\item[$(*)_2$] $\Omega_{\Lambda}= \bigcup_{k<\omega} \Omega_{\Lambda, k}$.

\item[$(*)_3$] $\Omega_{\Lambda,k}^+=\{ t ^{\frown} \langle 2n \rangle : t \in \Omega_{\Lambda,k}, n<\omega\}$.

\item[$(*)_4$] $\Omega_{\Lambda}^+=\bigcup_{k<\omega} \Omega_{\Lambda,k}^+$.

\item[$(*)_5$] $I_{\Lambda}= \{(t, f_t): t \in \Omega_{\Lambda}^+\}$, where $f_t \in \omega^\omega$ is defined as $f_t=t^{\frown} \langle 0: n<\omega \rangle$.

\item[$(*)_6$] $B_{\epsilon}=\{f \in \omega^{\omega} : t \in \Omega_{\Lambda}^+ \Rightarrow \neg (t \lhd f) \}$.
\end{enumerate}
We may note that each
$B_{\epsilon}$ is Borel, and  $\langle B_{\epsilon}: \epsilon <\omega_1  \rangle \in M.$
Furthermore $I_\Lambda \subseteq \mathbb{D}.$
\begin{claim}
$B_{\epsilon} \in I_{\mathbb D,\aleph_0}$.
\end{claim}
\begin{proof}
We have
\[
I_{\mathbb D,\aleph_0}=\{B \subseteq \omega^\omega: B \text{~contains a Borel set and~}\Vdash_{\mathbb D} \name{\eta}_{dom} \notin B  \}.
\]
Thus we have to show that $\eta_{dom} \notin B_{\epsilon}$.
For this, it is enough to show that $I_{\Lambda}$ is a maximal antichain, as then, it
will follow  that some $t \in \Omega_{\Lambda}^+$
is an initial segment of $\eta_{dom}$, and therefore, by the definition of $B_\epsilon$ we have $\eta_{dom} \notin B_{\epsilon}$, as requested.

First, we show  that $I_{\Lambda}$ is an antichain. Suppose
that $t \neq s$ are in  $\Omega_{\Lambda}^+$. We have to show that $(t, f_t)$ and $(s, f_s)$ are incompatible in $\mathbb D.$
 Suppose towards contradiction
that   $(t, f_t)$ and $(s, f_s)$ are compatible, in particular $t$ and $s$ are $\lhd$-comparable. Let us assume that
$t \lhd s$.
Let us write $t$ and $s$ as
\[
t= t_0 ^{\frown} \langle 2n_0+1\rangle ^{\frown} t_1 ^{\frown} \langle 2n_1+1 \rangle ^{\frown} \cdots ^{\frown} t_k ^{\frown} \langle 2n \rangle,
\]
and
\[
s=s_0 ^{\frown} \langle 2m_0+1\rangle ^{\frown} s_1 ^{\frown} \langle 2m_1+1 \rangle ^{\frown} \cdots ^{\frown} s_l ^{\frown} \langle 2m \rangle,
\]
where $n_i, m_j, n, m < \omega$ and $t_i, s_j \in \max(\Lambda)$.
The assumptions $t \lhd s$ and $t_i, s_j \in \max(\Lambda)$ imply that $k<l$, $s_i=t_i$ for $i \leq k$ and $m_i=n_i$ for $i<k$, hence we can write $s$ as
\[
s=t_0 ^{\frown} \langle 2n_0+1\rangle ^{\frown} t_1 ^{\frown} \langle 2n_1+1 \rangle ^{\frown} \cdots ^{\frown} t_k ^{\frown} \langle 2m_k+1 \rangle
^{\frown} \cdots ^{\frown} s_l ^{\frown} \langle 2m \rangle.
\]
Thus we must have $2n=2m_k+1$ which is impossible.
 Therefore, $I_{\Lambda}$
is an antichain.

Next we show that $I_{\Lambda}$ is a maximal antichain. Let
$p=(t_p, f_p) \in \mathbb D$. If there exists $t \in \Omega_{\Lambda}^+$
such that $t \unlhd t_p$, then $(t, f_t) \leq p$
and we are done. Therefore, we may assume that there is no such $t$.
Let
\[
\Omega^*=\{\langle \rangle\} \cup \{t: (\exists t_0,\cdots, t_{k-1} \in \max(\Lambda)) t=t_0 ^{\frown} \langle 2n_0+1\rangle ^{\frown} \cdots ^{\frown} t_{k-1}~ ^{\frown} \langle 2n_{k-1}+1\rangle \text{~and~} t \unlhd t_p\}.
\]
Then $\Omega^* \neq \emptyset$, and since $t_p$ has finite length,
there is an element $t$ of $\Omega^*$ of maximal length. Let
$s_1 \in \Lambda_{\epsilon}$ be such that $t ^{\frown} s_1 \unlhd t_p$
and $s_1$ is maximal. There are two cases.

{\bf Case I. $t ^{\frown} s_1 =t_p$.} Let $k$ be maximal such that
$$s_2 := s_1^{\frown} \langle f_p( \ell g(t)+\ell g(s_1)+i) : i<k \rangle \in \Lambda_{\epsilon}.$$
Note that by the construction,  $\Lambda_{\epsilon}$ has no infinite branches, and since $s_1 \in \Lambda_\epsilon,$ it follows that there is such maximal $k$.
It then follows from the choice of $s_2$ that $s_2 \in \max(\Lambda_\epsilon)$.
 Let
 $$s= t ^{\frown} s_1 ^{\frown} \langle f_p( \ell g(t)+\ell g(s_1)+i) : i<k \rangle=t_p ^{\frown} \langle f_p( \ell g(t)+\ell g(s_1)+i) : i<k \rangle.$$
then $t_p \unlhd s \in \Omega_{\Lambda}^+$ and $(s,f_s) \in I_\Lambda$
is compatible with $p$.

{\bf Case II. $t ^{\frown} s_1 \lhd t_p$.} First note that $s_1 \in \max(\Lambda_{\epsilon})$,
as otherwise, we will have $\Suc_{\Lambda_{\epsilon}}(t)=\{t^{\frown} \langle k \rangle : k<\omega\}$, and then for some
$k<\omega$, $t ^{\frown} s_1^{\frown} \langle k \rangle \unlhd t_p$ and
$s_1^{\frown} \langle k \rangle \in \Lambda_\epsilon,$ which contradicts the choice of $s_1$
as a maximal element of $\Lambda_\epsilon$ with  $t ^{\frown} s_1 \unlhd t_p$.
Now there are two possibilities:
\begin{itemize}
\item  $t_p( \ell g(t)+\ell g(s_1))$ is odd:
Then $t ^{\frown} s_1 ^{\frown} \langle t_p( \ell g(t)+\ell g(s_1)) \rangle \in \Omega^*$,
which contradicts the maximality of $t$.

\item  $t_p( \ell g(t)+\ell g(s_1))$ is even:  Then we have $t ^{\frown} s_1 ^{\frown} \langle t_p( \ell g(t)+\ell g(s_1)) \rangle  \in \Omega_\Lambda^+$
and  $t ^{\frown} s_1 ^{\frown} \langle t_p( \ell g(t)+\ell g(s_1)) \rangle \unlhd t_p$,
which contradicts our assumption that there is no element in $\Omega_\Lambda^+$
which is $\unlhd t_p$.
\end{itemize}
Thus case II cannot happen. The claim follows.
\end{proof}
We now turn to the main part of the lemma.
\begin{claim}
$\bigcup_{\epsilon<\omega_1}B_{\epsilon} \notin I_{\mathbb D,\aleph_0}$.
\end{claim}
\begin{proof} Suppose towards contradiction that $\bigcup_{\epsilon<\omega_1}B_{\epsilon} \in I_{\mathbb D,\aleph_0}$.
It then follows that  there exists a Borel set $B$ such that $\bigcup_{\epsilon<\omega_1}B_{\epsilon}  \subseteq B$
and $\Vdash_{\mathbb D}$``$\name{\eta}_{dom} \notin B$''.

For each condition $p \in \mathbb D$ set
\[
X(p)=\{g \in \omega^\omega: t_p \lhd g \text{~and~} f_p \leq g \}.
\]

By the definition of the ideal $I_{\mathbb D,\aleph_0}$, there exists a sequence
$\bar I=\langle I_n: n<\omega   \rangle$ of countable pre-dense subsets of $\mathbb D$
such that letting $I_n=\{p_{n,l} : l<\omega     \}$, we have
$B \subseteq \omega^\omega \setminus \big(\bigcap_{n<\omega} \bigcup_{l<\omega}X(p_{n,l})\big)$.

Let $\chi$ be a large enough regular cardinal and let $N$ be a countable elementary submodel of $L_{\chi}[\bar I, \langle \Lambda_{\epsilon} : \epsilon<\omega_1 \rangle]$
 such that $\bar I,\langle \Lambda_{\epsilon} : \epsilon<\omega_1 \rangle \in N$.
Let $\delta(*)=N\cap \omega_1^{L_{\chi}[\bar I, \langle \Lambda_{\epsilon} : \epsilon<\omega_1 \rangle]}$. Then
$\delta(*) < \omega_1$, and hence $\Lambda=\Lambda_{\delta(*)}$ is well-defined. Let $\bar I^*=\langle I^*_m : m<\omega \rangle$
enumerate all countable pre-dense subsets of $\mathbb D$ in $N$, and for every $m<\omega$,
set $I^*_m=\{ p_{m,l}^* : l<\omega \}$. Note that for every $n<\omega$, $I_n \in N$, hence
there exists some $j(n)<\omega$ such that $I_n=I^*_{j(n)}$.
We now prove the following.
\begin{itemize}
\item[$(*):$] There exists a sequence $\langle q_n; n<\omega      \rangle$
of conditions such that:

\begin{enumerate}
\item $q_n=(t_n,f_n) \in \mathbb D \cap N$.

\item $n=m+1 \Rightarrow q_m \leq q_n$.

\item If $n=m+1$, then there exists $l$ such that $p_{m,l}^* \leq q_n$.

\item $t_0=\langle \rangle$,
\item If $n>0$ then $t_n=s_n ^{\frown} \langle 2m_n \rangle$ for
some $s_n \in \Omega_{\Lambda}$ and $m_n<\omega$.
\end{enumerate}
\end{itemize}
We construct the sequence $\langle q_n; n<\omega      \rangle$ by induction on $n$.  For $n=0$ set $q_0=(\langle \rangle, \id_\omega)$,
where $\id_\omega$ is the identity function on $\omega.$ Now assume that $n=m+1$
and the condition $q_m=(t_m,f_m)$ is constructed, such that it satisfies the induction
hypothesis.
As $I^*_m=\{p_{m,l}^* : l<\omega\}$
is pre-dense, the set $E=\{p \in \mathbb D: (\exists l)(p_{m,l}^* \leq p)\} \in N$ is
open dense in $\mathbb{D}$. Let   $I=\{r_l : l<\omega\} \subseteq E$ be a maximal antichain in $\mathbb D$
above $q_m$. By
elementarity, we can assume that $I \in N$. Let $A=\{\tr(r_l) : l<\omega\}$. Then
by Lemma \ref{rank}(b) and the fact that $q_m, A \in N$, we have
\begin{itemize}
\item $\rk_{q_m,A}(t_{m})<\omega_1$,

\item $\rk_{q_m,A}(t_{m}) \in N$.
\end{itemize}
Thus $\rk_{q_m,A}(t_{m}) < \delta(*).$
Let $h_{\delta(*)} : \Lambda_{\delta(*)} \rightarrow \delta(*)+1$
witness $( \Lambda_{\delta(*)}, h_{\delta(*)}) \in Y_{\delta(*)}$,
 and let $\Lambda^*$
be the set of sequences $t \in \Lambda_{\delta(*)}$ satisfying
the following conditions:
\begin{itemize}
\item $q_m \leq (t_{m} ^{\frown} t, t_{m} ^{\frown} t ^{\frown} f_m \restriction [\ell g(t_{m} ^{\frown} t),\omega))$.

\item $\rk_{q_m, A}(t_{m} ^{\frown} t)< h_{\delta(*)}(t)$.
\end{itemize}
As $\delta(*)$ is a limit ordinal,
$\langle \rangle \in \Lambda^*$ and hence $\Lambda^* \neq \emptyset$. Let
$$\alpha_*=\min \{\rk_{q_m, A}(t_{m} ^{\frown} t) : t \in \Lambda^*\},$$
and choose $t_* \in \Lambda^*$ such that $\alpha_*=\rk_{q_m, A}(t_{m} ^{\frown} t_*)$.

We will show that $\alpha_*=0$. Assume towards a contradiction that
 $\alpha_*>0$. As $t_* \in \Lambda^*$, we have $\alpha_*=\rk_{q_m, A}(t_m^{\frown} t_*) < h_{\delta(*)}(t_*)$,
therefore by the definition of $h_{\delta(*)}$, for every $k$ large
enough, $\alpha_* \leq h_{\delta(*)}(t_* ^{\frown} \langle k \rangle)$.
By the definition
of the rank, the  set
$$U_1=\{k < \omega : \alpha_* \leq \rk_{q_m, A}(t_m ^{\frown}t_* ^{\frown} \langle k \rangle)\}$$
is finite. It is also clear that the set
$$U_2=\{k < \omega : k\leq f_m(\ell g(t_m^{\frown} t_*)) \}$$
is finite. It follows that the set $U=U_1 \cup U_2$ is finite, and hence, for every large enough $k$ we have $k \notin U$,
and
$\alpha_* <h_{\delta(*)}(t_* ^{\frown} \langle k \rangle)$. For such $k$,
\begin{itemize}
\item $\rk_{q_m, A}(t_m ^{\frown}t_* ^{\frown} \langle k \rangle) < \alpha_*$ (as $k \notin U_1$),
\item $\rk_{q_m, A}(t_m ^{\frown}t_* ^{\frown} \langle k \rangle) < \rk_{q_m, A}(t_m ^{\frown}t_*)$ (by the definition of the rank).
\end{itemize}
 Therefore, for every such $k$,
$t_* ^{\frown} \langle k \rangle \in \Lambda^*$ and $\rk_{q_m, A}(t_m ^{\frown}t_* ^{\frown} \langle k \rangle) < \alpha_*$,
which contradicts the minimality of $\alpha_*$.
Thus $\alpha_*=0$.

By the way we
defined the rank, we can find some $t' \in A$ such that:
\begin{itemize}
\item $q_m \leq (t',f_m \restriction [\ell g(t'),\omega))$,

\item $\ell g(t') \leq \ell g(t_{m} ^{\frown} t_*)$,

\item  $(t_{m} ^{\frown} t_*)(l) \leq t'(l)$,
for every $\ell g(t_m) \leq l <\ell g(t')$.
\end{itemize}
As $A=\{\tr(r_l) : l<\omega\}$ and $t' \in A,$ we have $t'=\tr(r_{l_*})$, for some
  $l_* < \omega$.
By induction on $i \geq \ell g(t')$ we choose some $k_i$ such that:
\begin{itemize}
\item $f_m(i) \leq k_i$,
\item $f_{r_{l_*}}(i) \leq k_i$,
\item  $s_{m,i}=t' \restriction [\ell g(t_m), \ell g(t')) ^{\frown} \langle k_j : \ell g(t') \leq j<i \rangle \in \Lambda_{\delta(*)}$.
\end{itemize}
Note that $t' \restriction [\ell g(t_m), \ell g(t')) \in \Lambda_{\delta(*)},$ and since
 $\Lambda_{\delta(*)}$ has no infinite branches,
it follows that there exists a maximal $i$ for which we can choose
$k_i$ as requested. Set
\begin{center}
$t''=t'$$^{\frown} \langle k_j : j<i+1 \rangle ^{\frown} \langle 2\big(f_{r_{l_*}}(\ell g(t')+i)+f_m(\ell g(t')+i)\big) \rangle.$
\end{center}
Let $q_{m+1}=(t'', t''$$^{\frown} f \restriction [\ell g(t''),\omega)),$
where $f(i)=max \{f_m(i),f_{r_{l_*}}(i)\}$ for every $i\in [\ell g(t''),\omega))$.
It's easy to see that $r_{l_*},q_m \leq q_{m+1}$
and that the condition $q_{m+1}$
satisfies the requirements (1)-(5) of $(*)$.
This completes the construction of the sequence
 $\langle q_n:  n<\omega      \rangle$.

 As the sequence  $\langle q_n=(t_n, f_n): n<\omega      \rangle$ is increasing, the sequence  $\langle t_n; n<\omega      \rangle$
is increasing too, and hence $f=\bigcup_{n<\omega}t_n$ is a
well-defined function.

We first show that $f \in \omega^\omega$. It suffices to show that $\dom(f)=\omega$. Thus let $k<\omega.$ As $N$ is an elementary submodel of $L_{\chi}[\bar I, \langle \Lambda_{\epsilon} : \epsilon<\omega_1 \rangle]$,
there is a pre-dense set $I\in N$ such that  for every $p \in I,$ $k <\ell g(t_p)$.
Let $m$ be such that $I=\{p_{m,l}^* : l<\omega\}$. By the way we defined $q_{m+1},$ for some $l,$
 $p_{m,l}^* \leq q_{m+1}$, hence $k< \ell g(t_{m+1})$. Thus $k \in \dom(t_{m+1}) \subseteq \dom(f),$ as requested.

Now we show that $f \notin B.$ It suffices to show that $f$  is $(N,\mathbb D)$-generic.  As the sequence $\langle q_n; n<\omega      \rangle$ is increasing,
it follows that $f \in X(q_n)$ for every $n<\omega$.
Let $I \in N$ be  a countable pre-dense subset of $\mathbb D$. Then for some $m$,
$I=I^*_m=\{p_{m,l}^* : l<\omega\}$.
Let $l$ be such that
$p_{m,l}^* \leq q_{m+1}$. It follows that
 $f \in X(p_{m,l}^*)$, and hence $t_{p_{m,l}^*} \lhd f$ and $f_{p_{m,l}^*} \leq f.$
 This implies  $f$  is $(N,\mathbb D)$-generic, as wanted.

We  now prove that $f \in B_{\delta(*)}$.
Let us recall that $B_{\delta(*)}=\{f \in \omega^{\omega} : t \in \Omega_{\Lambda_{\delta(*)}}^+ \Rightarrow \neg (t \lhd f) \}.$
Suppose towards contradiction that
$f \notin B_{\delta(*)}$.
Then for some $t \in \Omega_{\Lambda_{\delta(*)}}^+$,
we have $t \lhd f.$ Let $n$ be large enough such that
$t \lhd t_n \lhd f.$ Then $(t, f_t), (t_n, f_{t_n}) \in I_{\Lambda_{\delta(*)}}$
are compatible, as witnessed by $(t_n, f)$, which contradicts
the fact that $I_{\Lambda_{\delta(*)}}$ is an antichain.

But then $f \in B_{\delta(*)}$ and $B_{\delta(*)} \subseteq B$, and hence $f \in B$ which is a contradiction.
The claim follows.
\end{proof}
This completes the proof of Lemma \ref{positive}.
\end{proof}
We now extend the above result to include the class of all Suslin ccc forcing notions which add a Hechler real.
\begin{lemma} ($\ZF+DC$)
\label{positive-g}
Let $M$ be an inner model of $\ZFC$ with
$\aleph_1^M=\aleph_1$. Let $\mathbb Q \in M$ be a Suslin ccc forcing notion
which adds a  Hechler real,
and let $I=I_{\mathbb Q, \aleph_0}$.
 Then the union of all Borel $I$-small sets coded in the model $M$ is $I$-positive.
\end{lemma}
\begin{proof} Let $\name \eta$ be the canonical name for a real added by $\mathbb Q$ and let $f$ be a Borel function such that $\Vdash_{\mathbb Q}$``$f(\name{\eta})=\name{\eta}_{dom}$''. We may assume that $\name \eta, f \in M.$ Let also  $\langle B_{\alpha} : \alpha<\omega_1 \rangle$ be the sequence constructed in the
proof of Lemma \ref{positive}.
Then the sequence $\langle f^{-1}(B_{\alpha}) : \alpha<\omega_1 \rangle$ is in $M$.
\begin{claim}
For every $\alpha < \omega_1, f^{-1}(B_{\alpha}) \in I_{\mathbb Q,\aleph_0}.$
\end{claim}
\begin{proof}
Let $\alpha < \omega_1.$ Then $\Vdash_{\mathbb Q}$``$\name{\eta}_{dom} \notin B_{\alpha}$'',
and therefore $\Vdash_{\mathbb Q}$``$\name{\eta} \notin f^{-1}(B_{\alpha})$''. Thus
  $f^{-1}(B_{\alpha}) \in I_{\mathbb Q,\aleph_0}$, as required.
 \end{proof}
 \begin{claim}
 $\bigcup_{\alpha < \omega_1} f^{-1}(B_{\alpha}) \notin I_{\mathbb Q,\aleph_0}.$
\end{claim}
\begin{proof}
Let $\chi$ be large enough regular and let $N \prec (H(\chi), \in)$ be a countable elementary submodel of $H(\chi)$ containing all the relevant
objects. It suffices to find $K \subseteq \mathbb Q \cap N$ such that
$K$ is $(N,\mathbb Q)$-generic and $\name{\eta}[K] \in \bigcup_{\alpha<\omega_1}f^{-1}(B_{\alpha})$.

Let $RO(\mathbb D)$ be the Boolean completion of $\mathbb D$ and let $\pi: \mathbb Q \rightarrow RO(\mathbb D)$ be a projection, which exists by our assumption on $\mathbb Q.$
By Lemma \ref{positive}, there exists  $G\subseteq RO(\mathbb D) \cap N$ which is $(N,RO(\mathbb D))$-generic, such
that $\name{\eta}_{dom}[G] \in \bigcup_{\alpha<\omega_1}B_{\alpha}$.
Let $H \subseteq \mathbb Q / G$ be $\mathbb Q / G$-generic over $N[G],$ where
$\mathbb Q / G=\{q \in \mathbb Q: \pi(q) \in G  \}$.
Then $K=G \ast H$ is  $(N,\mathbb Q)$-generic and
clearly
$\name{\eta}[K] \in \bigcup_{\alpha<\omega_1}f^{-1}(B_{\alpha})$.
\end{proof}
The lemma follows.
\end{proof}
\section{A measurable cardinal from regularity properties
and $\DC_{\omega_1}$}
\label{regularity}

In this section we prove a general criterion for the existence
of an inner model for a measurable cardinal under the assumptions
``$\ZF+\DC_{\omega_1}+$ all sets of reals have certain regularity properties''.

Let us recall that an ideal $J$ on a set $X$ is $\kappa$-saturated, if $P(X)/
J,$ considered as a forcing notion, satisfies the $\kappa$-cc.  The following is well-known.
\begin{lemma} (see \cite{kunen})
\label{l-2} Assume
there exists a non-trivial $\kappa$-complete $\kappa$-saturated ideal. Then there exists an inner model with a measurable cardinal.
\end{lemma}
We now prove the following general result.
\begin{lemma}  ($\ZF+\DC$)
\label{l-3}
Assume  $\lambda$ is an uncountable cardinal  and the following conditions hold:
\begin{enumerate}
\item $\mathbb Q$ is a  Suslin ccc forcing notion,

\item $I$ is a $\sigma$-complete ideal on the reals extending $I_{\mathbb Q,\aleph_0}$,

\item $\langle B_{\alpha} : \alpha<\lambda \rangle$ is a sequence of sets from $I_{\mathbb{Q}, \aleph_0}$
such that $\bigcup_{\alpha < \lambda} B_{\alpha}  \notin I$,

\item There is no sequence $\langle B_{\alpha}^* : \alpha<\aleph_1 \rangle$ of
$I$-positive sets such that $B_{\alpha}^* \cap B_{\beta}^* \in I$
for every $\alpha \neq \beta <\aleph_1$.
\end{enumerate}
Then there exists an inner model of $ZFC$ with a measurable cardinal.
\end{lemma}
\begin{remark}
 \label{rem-1} The condition in clause (4) follows from the assumption $\mathcal{P}(\omega^{\omega})/I \models ccc$, and these are equivalent under $\text{AC}_{\omega_1}$.
\end{remark}
\begin{proof}
For every $\alpha<\lambda$, let $B'_{\alpha}=B_{\alpha} \setminus \bigcup_{\beta<\alpha}B_{\beta}$.
By clause (3), $\bigcup_{\alpha < \lambda} B'_{\alpha}=\bigcup_{\alpha < \lambda} B_{\alpha}  \notin I.$ Let $U \subseteq \lambda$
be of minimal size such that $\bigcup_{\alpha \in U} B'_{\alpha} \notin I.$
Let also  $(\xi_{\alpha} : \alpha<|U|)$ enumerate $U$ and for each $\alpha<|U|$ set
$B''_{\alpha}=B'_{\xi_{\alpha}}$. Then $\langle B''_{\alpha}: \alpha<|U| \rangle$
is a sequence of pairwise disjoint sets whose union is $I$-positive.
Let $J$ be the ideal on $|U|$ defined as:
\[
X \in J \Leftrightarrow \bigcup_{\alpha \in X} B''_\alpha \in I.
\]
\begin{claim}
Let $\kappa=\aleph_1^V.$ Then
 $L[J] \models$``$J \cap L[J]$ is a $\kappa$-complete $\kappa$-saturated ideal
on $|U|$''.
\end{claim}
\begin{proof}
We first show that $J$ is $\kappa$-complete in $V$. Suppose towards a contradiction that $\langle X_n : n<\omega \rangle \in V$
is such that $X_n \in J$ for $n<\omega$, but $X=\bigcup_{n<\omega}X_n \notin J$.
For $n<\omega$ set $A_n=\bigcup_{\alpha \in X_n} B_{\alpha}''$. Then $\langle A_n : n<\omega \rangle \in V$,
and by the definition of the ideal $J$,
\begin{itemize}
\item $A_n \in I$, for $n<\omega,$
\item $A=\bigcup_{n<\omega}A_n \notin I$.
\end{itemize}
This contradicts clause (2). It immediately follows that
$L[J] \models$``$J \cap L[J]$ is a $\kappa$-complete ideal on $|U|$''.

We now show that $L[J]\models$``$\mathcal{P}(|U|)/(J \cap L[J])$ is $\kappa$-cc''.
Suppose not. Thus  there exists a sequence $\langle A_{\alpha} : \alpha<\kappa \rangle \in L[J]$
of $J$-positive sets such that for every $\alpha<\beta < \kappa$, $A_{\alpha} \cap A_{\beta} \in J$.
As $J$ is $\kappa$-complete in $L[J]$, we may assume without loss of generality that
$A_{\alpha} \cap A_{\beta}=\emptyset,$ for every $\alpha<\beta < \kappa$.

Work in $V$. For each $\alpha < \kappa$ set $B^*_\alpha=\bigcup_{\xi \in A_\alpha}B''_\xi.$
Then $\langle B^*_\alpha: \alpha < \kappa \rangle \in V$ satisfies the following conditions:
\begin{itemize}
\item $B^*_\alpha \notin I$ for every $\alpha < \kappa$ (as $A_\alpha \notin J$),

\item $B^*_\alpha \cap B^*_\beta=\emptyset$, for $\alpha < \beta < \kappa.$
\end{itemize}
This contradicts clause (4).
The claim follows.
\end{proof}
Thus by Lemma \ref{l-2},  there exists  an inner model of $\ZFC$ with a measurable cardinal.
\end{proof}
Assume $\ZF+\DC$, and let $\mathbb Q$ be a Suslin ccc forcing notion which adds a Hechler real. Then the ideal  $I_{\mathbb Q, \aleph_0}$
 is a $\sigma$-complete ideal and by Lemma \ref{positive-g},
it  is not $\aleph_2$-complete. We now show that the forcing notion $\mathcal{P}(\omega^{\omega})/I_{\mathbb Q,\aleph_0}$
satisfies the ccc. In order to do this,
we use $\DC_{\omega_1}$
and some regularity properties.

The following definition is of interest in the absence of choice.

\begin{definition}
Assume $\mathbb Q$ is a forcing notion, $B$ is a Boolean algebra and $I \subseteq B$ is an ideal.
\begin{enumerate}
\item We say that  $\mathbb Q$
satisfies the strong chain condition (scc), if there is no uncountable%
\footnote{So it may be non well-orderable in the absence of choice. %
} set $\{X_s : s\in S\} \subseteq \mathcal{P}(\mathbb Q)$ such that:
\begin{enumerate}
\item $X_s \neq \emptyset$ for each $s\in S$,
\item for every $s\neq t$ in  $S$,
if $p\in X_s$ and $q\in X_t$, then $p$ and $q$ are incompatible.
\end{enumerate}
\item We say that the pair $(B,I) $ satisfies the weak strong chain condition (scc$^-$)
if there is no uncountable collection $\{X_s : s\in S\} \subseteq \mathcal{P}(B)$
of non-empty subsets of $B$ such that:
\begin{enumerate}
\item  $X_s \cap I=\emptyset,$ for each $s\in S$,
\item for every $s\neq t$ in $S$, if  $b_s \in X_s$ and  $b_t \in X_t$, then $b_s \wedge b_t \in I$.
\end{enumerate}
\item We say that the pair $(B,I)$
satisfies the weak countable chain condition (ccc$^-$) if there is no uncountable collection
$\{b_s : s\in S\} \subseteq B$ of $I$-positive elements of $B$ such that $b_s \cap b_t \in I$,
for every $s \neq t$ in $S$.
\end{enumerate}
\end{definition}
\begin{remark}
Given an infinite cardinal $\kappa,$ we can similarly define the notions of $\kappa$-strong chain condition,
$\kappa$-weak strong chain condition and $\kappa$-weak  chain condition.
\end{remark}
We have the following easy lemma.
\begin{lemma}
\label{l-4} Assume $\mathbb Q$ is a forcing notion and $\DC_{\omega_1}$ holds.
Then $\mathbb{Q}$ satisfies the  strong chain condition if and only if it satisfies the ccc.
\end{lemma}
Let $\Borel(\omega^{\omega})$
denote the collection of all Borel subsets of $\omega^{\omega}$.
\begin{lemma}
\label{l-5}
($\ZF$)
Let $\mathbb Q$ be a Suslin  forcing
notion, and suppose that it satisfies the strong chain condition. Then:
\begin{enumerate}
\item[(a)]  $(\Borel(\omega^{\omega}),I_{\mathbb Q,\aleph_0})$ satisfies the weak countable chain condition.
\item[(b)] $(\Borel(\omega^{\omega}),I_{\mathbb Q,\aleph_0})$ satisfies  the weak strong chain condition.
\end{enumerate}
\end{lemma}
\begin{proof}
For notational simplicity set $\mathcal{B}=\Borel(\omega^{\omega})$.

(a).
Let $\name\eta$ be the canonical $\mathbb{Q}$-name for a real.
Suppose that $\{B_s : s\in S\} \subseteq \mathcal{B}$ is a collection of $I_{\mathbb Q,\aleph_0}-$positive
Borel sets such that for $s\neq t $ in $S$,   $B_s \cap B_t \in I_{\mathbb Q,\aleph_0}$.
For every $s\in S$, let
$$X_s=\{p\in \mathbb Q : p\Vdash \name{\eta} \in B_s\}.$$
As each $B_s$ is $I_{\mathbb Q,\aleph_0}-$positive, $X_s \neq \emptyset$. Furthermore if $s\neq t$ are in $S$, then since $B_s \cap B_t \in I_{\mathbb Q,\aleph_0}$,
for $p\in X_s$ and $q\in X_t$,  $p$ and $q$ are incompatible.
By our assumption, $\mathbb Q$ satisfies the strong chain condition, and hence   $S$ is countable. It follows that
$(\mathcal{B},I_{\mathbb Q,\aleph_0})$ satisfies the weak countable chain condition.

(b). Suppose that $\{X_s : s\in S\} \subseteq \mathcal{P}(\mathcal{B})$
is a collection of non-empty subsets of $\mathcal B$ such that for each $s \in S,~ X_s \cap I_{\mathbb Q,\aleph_0}=\emptyset$,
 and for  $s\neq t$ in $S$,  if $B_s\in X_s$ and $B_t \in X_t$, then $B_s \cap B_t \in I_{\mathbb Q,\aleph_0}$.
For each $s\in S$ let
\[
P_s=\{ p\in \mathbb{Q}: (\exists B \in X_s) p \Vdash\text{``}\name\eta \in B \text{''}        \}.
\]
Then each $P_s$ is a non-empty set and as in the proof of (a),  if $s\neq t$ are in $S$,
$p\in P_s$ and $q\in P_t$, then  $p$ and $q$ are incompatible. Thus, by our assumption, $S$ is countable.
The result follows.
\end{proof}

\begin{lemma}
\label{l-6}
($\ZF+\DC_{\omega_1}$) Assume  $\mathbb Q$ is
a Suslin ccc forcing notion. Then
\begin{enumerate}
\item[(a)] $\Borel(\omega^{\omega})/I_{\mathbb Q,\aleph_0}$ is ccc.
 \item[(b)] Assume all sets of reals are $I_{\mathbb Q,\aleph_0}$-measurable. Then $\mathcal{P}(\omega^{\omega})/I_{\mathbb Q,\aleph_0}$ is ccc.
 \end{enumerate}
\end{lemma}
\begin{proof}
(a) By Lemma \ref{l-4}, $\mathbb Q$ satisfies the strong chain condition. By Lemma \ref{l-5}(a),
$(\Borel(\omega^{\omega}), I_{\mathbb Q,\aleph_0})$ satisfies the weak countable chain condition. Now suppose by the way of contradiction that
 $\Borel(\omega^{\omega})/I_{\mathbb Q,\aleph_0}$ does not satisfy the ccc, and let
 $A=\{\mathcal X_s : s\in S\} \subseteq \Borel(\omega^{\omega})/I_{\mathbb Q,\aleph_0}$
be an uncountable antichain. Note that each $\mathcal X_s$ is an equivalence class.
By $\DC_{\omega_1}$  $S$
has a subset $\{s_\xi: \xi < \omega_1\}$ of size $\aleph_1$.
By another application of $\DC_{\omega_1}$, we can choose representatives $X_\xi \in \mathcal X_{s_\xi}$.
Then $\{ X_\xi: \xi < \omega_1   \} \subseteq  \Borel(\omega^{\omega})$ satisfies the following:
\begin{itemize}
\item Each $X_\xi$ is $I_{\mathbb Q,\aleph_0}$-positive,

\item For $\xi < \zeta < \omega_1,$ $X_\xi \cap X_\zeta \in I_{\mathbb Q,\aleph_0}.$
\end{itemize}
This contradicts the fact that the pair $(\Borel(\omega^{\omega}), I_{\mathbb Q,\aleph_0})$
satisfies the weak countable chain condition.

(b). By the way of contradiction suppose that
$\{X_s : s\in S\}$ is an uncountable collection of $I_{\mathbb Q,\aleph_0}$-positive
sets such that for $s\neq t$ in $S, X_s \cap X_t \in I_{\mathbb Q,\aleph_0}$.
 For
each $s\in S$, let
\begin{center}
$P_s=\{B\subseteq \omega^{\omega} : B$ is a
Borel set such that $B=X_s$ mod $I_{\mathbb Q,\aleph_0}\}$.
\end{center} By our assumption, each $P_s$
is non-empty. By $\DC_{\omega_1}$, there
exists an uncountable subset $\{s_\xi: \xi < \omega_1   \}$
of $S$, and by another application of $\DC_{\omega_1}$, we can choose the representatives $B_\xi \in P_{s_\xi},$
for $\xi < \omega_1$. Then the collection $\{B_\xi: \xi < \omega_1    \}$
witnesses that $\Borel(\omega^{\omega})/I_{\mathbb Q, \aleph_0}$
does not satisfy ccc, which contradicts (a).

\end{proof}
The following is an immediate corollary of the above results.
\begin{corollary} ($\ZF+\DC_{\omega_1}$)
\label{c-1} Assume $\mathbb Q$ is a Suslin ccc forcing notion, $\lambda$ is an infinite cardinal and the following conditions hold:
\begin{enumerate}
\item All sets of reals are $I_{\mathbb Q,\aleph_0}-$measurable,

\item There exists a sequence $\langle B_{\alpha} : \alpha<\lambda \rangle$ of sets from
$I_{\mathbb{Q},,\aleph_0}$, such that $\bigcup_{\alpha < \lambda}B_{\alpha} \notin I_{\mathbb Q, \aleph_0}$.
\end{enumerate}
Then there is an inner model of $\ZFC$ with a measurable cardinal.
\end{corollary}
\begin{proof}
Let $I=I_{\mathbb Q, \aleph_0}$. By Lemma \ref{l-3} and Remark \ref{rem-1}, it is enough to show that
$\mathcal{P}(\omega^{\omega})/I$ satisfies the ccc.
This follows from Lemma \ref{l-6}(b).
\end{proof}

\begin{theorem} ($\ZF+\DC_{\omega_1}$)
\label{c-2}
 Let $\mathbb Q$ be a Suslin ccc forcing
notion which adds a Hechler real.
Suppose every set of reals is $I_{\mathbb Q,\aleph_0}$-measurable.
Then there is an inner model of $\ZFC$ with a measurable cardinal.
\end{theorem}
\begin{proof}
By Lemma \ref{positive-g} and Corollary \ref{c-1}.
\end{proof}

\section{A class of Suslin ccc forcing notions adding a Hechler real}
\label{prel}
In this section we prove Theorem \ref{t2}, by showing that the
 forcing notions $\mathbb{Q}_{\bold n}^1$ from  \cite{HwSh:1067}  add a Hechler real.
\begin{remark}
In \cite{HwSh:1067}, some other classes of forcing notions were also introduced, but as they are not relevant to our work, we do not discuss them here.
\end{remark}
Let us recall some definitions and facts from \cite{HwSh:1067}.
\begin{definition}
 A norm on a set $A$ is a function
 \[
 \nor: \mathcal{P}(A)\setminus \{\emptyset\} \rightarrow [0, \infty)
 \]
such that if  $X \subseteq Y$, then $\nor(X) \leq \nor(Y)$.
\end{definition}
In order to define the forcing notions  $\mathbb{Q}_{\bold n}^1$, we first need to define the corresponding parameters $\bold n$.

\begin{definition}
\label{def1}
A \emph{nice parameter} is a tuple
$\bold{n}=(T, \nor, \bar{\lambda}, \bar{\mu})$
such that:
\begin{enumerate}
\item $T$ is a subtree of $\omega^{<\omega}$,

\item $\bar{\mu}=(\mu_{t} :t \in T)$ is a sequence of non-negative
real numbers,

\item $\bar{\lambda}=(\lambda_{t} : t \in T)$ is a sequence of
pairwise distinct non-zero natural numbers such that for each $s, t \in T$:
\begin{enumerate}
\item  $\ell g(t) <  \mu_{t}  < \lambda_{t}=|\Suc_T(t)|$,\footnote{It follows that $T\cap \omega^n$
is finite and non-empty for every $n>0$.}

\item If $\ell g(s)=\ell g(t)$ and $s <_{\lex} t$ then $\lambda_{s} < \lambda_{t}$,

\item If $\ell g(s)< \ell g(t)$,  then $ \lambda_{s} <  \lambda_{t}$,
\end{enumerate}

\item For $t \in T$, $\nor_{t}$ is a norm on $\Suc_T(t)$
such that:
\begin{enumerate}

\item $(\ell g(t)+1)^2 \leq \mu_{t}\leq \nor_{t}(\Suc_T(t))$,

\item $\lambda_{<t} < \mu_{t}$, where $\lambda_{<t}= \prod_{ \lambda_{s}<\lambda_{t}}\lambda_{s},$

\item (Co-Bigness) Suppose $r>0$, $i(*)\leq \mu_{t}$ and for every $i< i(*)$,  $a_i \subseteq \Suc_{T}(t)$
and $r+\frac{1}{\mu_{t}}\leq \nor_{t}(a_i)$.
 Then $r\leq \nor_{t}(\bigcap_{i<i(*)}a_i)$.

\item If $1\leq \nor_{t}(a)$ then $\frac{1}{2}<\frac{|a|}{|\Suc_{T}(t)|}$.

\item If $r+ \mu_{t}\leq \nor_{t}(a)$ and $s \in a$, then $r\leq \nor_{t}(a\setminus \{s\})$.
\end{enumerate}
\end{enumerate}
\end{definition}
\begin{notation}
Given a nice parameter $\bold n,$ we denote it as $\bold n=(T_{\bold n}, \nor_{\bold n}, \bar{\lambda}_{\bold n}, \bar{\mu}_{\bold n})$,
where $\bar{\lambda}_{\bold n}= \langle \lambda^{\bold n}_t: t \in T_{\bold n}      \rangle$ and
 $\bar{\mu}_{\bold n}= \langle \mu^{\bold n}_t: t \in T_{\bold n}      \rangle$.
 Furthermore, we denote $(\nor_{\bold n})_t$ as $\nor^{\bold n}_t$.
\end{notation}
We are now ready to define the forcing notions $\mathbb{Q}^1_{\bold n},$ where $\bold n$ is a nice parameter.
\begin{definition}
\label{forcing}
Suppose $\bold n$ is a nice parameter. The forcing notion $\mathbb{Q}^1_{\bold n}$
is defined as follows.
\begin{enumerate}
\item $p\in \mathbb{Q}_{\bold n}^1$ iff for some $\tr(p) \in T_{\bold n}$
we have:
\begin{enumerate}
\item $p=(\tr(p), T_p)$, where  $T_p$ is a subtree of $T_{\bold n}$ with trunk $\tr(p)$,

\item For $\eta \in \lim(T_p)$,
$$\lim(nor_{\eta \restriction l}(\Suc_{T_p}(\eta \restriction l)) : \ell g(\tr(p)) \leq l<\omega)=\infty,$$

\item $2-\frac{1}{\mu_{\tr(p)}} \leq \nor(p)$, where
\[
\nor(p)=\sup \{a >0 : t \in T_p \Rightarrow a\leq \nor_{t}(\Suc_{T_p}(t))\},\footnote{Note that $\nor(p)=\inf \{\nor_{t}(\Suc_{T_p}(t)) : t \in T_p\}.$}
\]
\item For every $n<\omega$, there exists $k^p(n)> \ell g(\tr(p))$ such
that
\begin{center}
 $t \in T_p$ and $\ell g(t) \geq k^p(n)  \Rightarrow n\leq \nor_{t}(\Suc_{T_p}(t))$.
\end{center}
\end{enumerate}
\item Suppose $p, q \in \mathbb{Q}_{\bold n}^1$. Then $p \leq q$ iff $T_q \subseteq T_p$
\end{enumerate}
\end{definition}
To each forcing notion $ \mathbb{Q}_{\bold n}^1$ we can assign a canonical name for a real.
\begin{definition}
Suppose $\bold n$ is a nice parameter. Let $\name{\eta}_{\bold n}^{1}$
be the $\mathbb{Q}_{\bold n}^1-$name
$$\name{\eta}_{\bold n}^{1}=\bigcup\{\tr(p) : p\in \dot{G}_{\mathbb{Q}_{\bold n}^1}\},$$
where $ \dot{G}_{\mathbb{Q}_{\bold n}^1}$ is the canonical $\mathbb{Q}_{\bold n}^1$-name for the generic filter.
\end{definition}
\begin{remark}
If $G$ is a $\mathbb{Q}_{\bold n}^1$-generic filter over $V$, then $\name{\eta}_{\bold n}^{1}[G]$
is a generic real added by $G$, furthermore, we can recover $G$ from $\name{\eta}_{\bold n}^{1}[G]$ by
\[
G=\{p \in \mathbb{Q}_{\bold n}^1: \tr(p) \lhd  \name{\eta}_{\bold n}^{1}[G] \}.
\]
\end{remark}
We  now state some of the basic properties and results on
$\mathbb{Q}_{\bold n}^1$.

\begin{theorem}
Suppose $\bold n$ is a nice parameter. Then:
\begin{enumerate}
\item[(a)] $\mathbb{Q}_{\bold n}^1$
is a Suslin ccc forcing notion.

\item[(b)] Forcing with $\mathbb{Q}_{\bold n}^1$ adds a Cohen real.

\item[(c)] The following is consistent relative to $\ZFC$:
\begin{enumerate}
\item[(1)] $\ZF$,

\item[(2)] Every set of reals is $I_{\mathbb{Q}_{\bold n}^1,\aleph_1}$-measurable,

\item[(3)] There exists an $\omega_1-$sequence of distinct reals.
\end{enumerate}
\end{enumerate}
\end{theorem}
\begin{proof}
See \cite{HwSh:1067}.
\end{proof}
We  now turn to the proof of Theorem \ref{t2}, and show that for each nice parameter $\bold n$, the forcing notion $\mathbb{Q}_{\bold n}^1$ adds a Hechler real.
In order to do that, we first prove the weaker result that $\mathbb{Q}_{\bold n}^1$
adds a dominating real.

\begin{lemma}
\label{l-7}
Suppose $\bold n$ is a nice parameter. Then forcing with $\mathbb{Q}_{\bold n}^1$ adds a dominating
real.
\end{lemma}
\begin{proof}
For every $t \in T_{\bold n}$ and $k\leq \ell g(t)$,
let  $w_{t,k} \subseteq \Suc_{T_{\bold n}}(t)$  be such that:
\begin{itemize}
\item $\nor_{t}(\Suc_{T_{\bold n}}(t) \setminus w_{t, k})=k+1$,
\item $|w_{t, k}|$ is
minimal,
\item  If $k+1 \leq \ell g(t),$ then  $w_{t,k+1} \subseteq w_{t,k}$.
\end{itemize}
 Then we have the following.
 \begin{claim}
 \label{cla1}
 Let $t$ and $k$ be as above.
\begin{enumerate}
\item[(a)]  If $u\subseteq \Suc_{T_{\bold n}}(t)$ and
$k+2 \leq \nor_{t}(u)$, then $u\cap w_{t, k} \neq \emptyset$.

\item[(b)] If $u\subseteq \Suc_{T_{\bold n}}(t)$, $l<k$ and $l+1 \leq \nor_{t}(u)$,
then letting $v=u\setminus w_{t, k}$, we have:
\begin{enumerate}
\item[(1)] $v\subseteq u$ and $v\cap w_{t, k}=\emptyset$.

\item[(2)]  $l\leq \nor_{t}(v)$
and $v\neq \emptyset$.

\item[(3)] If $\nor_{t}(u)>2$,
then $\min\{k, \nor_{t}(u)-1\} \leq \nor_{t}(v)$.
\end{enumerate}
\end{enumerate}
\end{claim}
\begin{proof}
(a). Otherwise,  $u\subseteq Suc_{T_{\bold n}} \setminus w_{t, k}$,
and hence $\nor_{t}(\Suc_{T_{\bold n}} \setminus w_{t, k}) \geq \nor_{t}(u) \geq k+2$,
which is impossible.

(b). Clause (1) is clear.
For clause (2), note that
$v=u \cap (\Suc_{T_{\bold n}}(t) \setminus w_{t, k})$,
hence by the co-bigness property \ref{def1}(4)(c),
$\nor_{t}(v)=\nor_{t}(u \cap (\Suc_{T_{\bold n}}(t) \setminus w_{t, k})) \geq l.$
For clause (3), let $j=\min\{k-1, \nor_{t}(u)-2\}=\min\{\nor_{t}(\Suc_{T_{\bold n}}(t) \setminus w_{t, k})-2, \nor_{t}(u)-2\}$.
By the co-bigness property,
$\nor_{t}(v)=\nor_{t}(u \cap (\Suc_{T_{\bold n}}(t) \setminus w_{t, k})) \geq j+1$,
from which the result follows.
\end{proof}
By induction on $n$ we define a sequence $\langle \name\tau_n: n<\omega      \rangle$ of $\mathbb{Q}_{\bold n}^1$-names for a member of $\omega \cup \{\omega\}$
as follows:
\begin{enumerate}
\item If $n=0$, then $\name\tau_0=\check{0}$,

\item If $n=m+1$ and $\name\tau_m=\check\omega$, then $\name\tau_n=\check\omega$ as well.
Otherwise, we let $\name\tau_n=\check{j}$ where $j$ is the
minimal natural number such that $\Vdash_{\mathbb{Q}_{\bold n}^1}$``$\name\tau_m < \check{j}$ and $\name\eta^1_{\bold n} \restriction j+1 \in \check{w}_{\name\eta^1_{\bold n} \restriction j, n}$'',
if such a $j$ exists. Otherwise, we let $\name\tau_n=\check\omega$.
\end{enumerate}
\begin{claim}
\label{cla2}
 $\Vdash_{\mathbb{Q}_{\bold n}^1}$``$\name\tau_n < \check{\omega}$, for every $n<\omega$.
\end{claim}
\begin{proof}
We prove the claim by induction on $n$. For $n=0$ the claim is obvious by the choice of $\name\tau_0$. Now suppose that
$n=m+1$ and the claim holds for $m$. Let $p\in \mathbb{Q}_{\bold n}^1$ be an arbitrary condition. We find an extension $q$ of
$p$, forcing $\name\tau_n < \check{\omega}$. By extending $p$ if necessary, we may assume that:
\begin{itemize}
\item $p$ decides $\langle \name\tau_i: i < n \rangle$, say it forces ``$\name\tau_i=j_i$'' for every $i<n$,

\item $j_m+m+1 < \ell g(\tr(p))$,

\item $n+2 < \nor_{t}(\Suc_{T_p}(t))$ for every $\tr(p) \unlhd t \in T_p$.
\end{itemize}
By Claim \ref{cla1}, $w_{\tr(p), n} \cap \Suc_{T_p}(\tr(p)) \neq \emptyset.$ Let $t \in w_{\tr(p), n} \cap \Suc_{T_p}(\tr(p))$
and let $q=(t, T_p^{[t]}),$ where $T_p^{[t]}=\{s \in T_p: s \unlhd t \text{~or~} t \unlhd s  \}$.
Then
\begin{center}
 $q \Vdash$``$\name\eta^1_{\bold n} \restriction (\ell g(\tr(p))+1)=t \in \check{w}_{\tr(p), n}= \check{w}_{\name\eta^1_{\bold n} \restriction \ell g(\tr(p)), n}$'',
 \end{center}
and hence  $q \Vdash$``$\name\tau_n \leq \ell g(\tr(p))$'', as required.
\end{proof}
Let $\name\eta$ be a name for a real  such that
\[
\Vdash_{\mathbb{Q}_{\bold n}^1}\text{``} (\forall n<\omega) \name\eta(n)=\name\tau_n.
\]
\begin{claim}
\label{cla3}
Suppose $h\in \omega^{\omega}$ and $p\in \mathbb{Q}_{\bold n}^1$.
Then there exists a condition $q \geq p$ such that $q\Vdash$``$h(n) \leq \name\eta(n)$
for every large enough $n$''.
\end{claim}
\begin{proof} Without loss of generality, we may assume that the function $h$ is increasing. By extending $p$ is necessary, we may
also assume that
$\ell g(\tr(p)) > 2$ and for every $t \in T_p$
with $\tr(p) \unlhd t$ we have $2 < \nor_{t}(\Suc_{T_p}(t))$.

Let $m$ be maximal such that $p$ decides $\name\eta \upharpoonright m+1$, and let $\langle j_i: i < m+1  \rangle$
be such that $p \Vdash$``$\name\eta \restriction m+1=\langle j_i: i < m+1  \rangle$''.
Let also $\langle n_i: i<\omega   \rangle$ be an increasing sequence of natural numbers such that;
\begin{itemize}
\item $n_0=\ell g(\tr(p)),$

\item $h(n_i) < n_{i+1},$ for every $i<\omega.$
\end{itemize}
Let
\[
T=\{t \in T_p: (\forall i<\omega)(\forall l\in [n_i,n_{i+1})) \big[ l < \ell g(t) \Rightarrow  t \restriction (l+1) \notin w_{t \restriction l,i}\big]\}.
\]
$T$ is obviously  downwards closed and by the co-bigness property, $T$ is a perfect tree.
Let $t \in T$ be such that $2< \nor_{s}(\Suc_{T}(s))$ for
every $t \unlhd s \in T$
and set $q=(t, T^{[t]})$. Then $q \in \mathbb{Q}_{\bold n}^1$ and $q \geq p$. We show that
\[
q \Vdash\text{``}  m < n \Rightarrow h(n) \leq \name\eta(n)=\name\tau_n        \text{''.}
\]
Suppose not. Thus we can find $r \geq q$ and $n> m$ such that
\begin{center}
$(*)_1$ \quad\quad $r \Vdash$``$\name\tau_n  < h(n)$''.
\end{center} By extending $r$, we may assume that:
\begin{itemize}
\item $r$ decides $\name\eta \upharpoonright n+1$, say it forces ``$\name\eta \upharpoonright n+1= \langle j_i: i< n+1 \rangle$'',

\item $j_n < \ell g(\tr(r))$.
\end{itemize}
Let $i$ be such that $n_i \leq j_n < n_{i+1}$. Set $s_1= \tr(r) \restriction j_n$ and $s_2=\tr(r) \restriction (j_n+1).$
Then $s_2 \in \Suc_{T_q}(s_1)$, and hence by the definition of the tree $T$ we have $s_2 \notin w_{s_1, i}$.
On the other hand, $r \Vdash$``$\name\tau_n=j_n < \ell g(\tr(r))$'', and hence
$s_2 \in w_{s_1, n}$. As the sequence $\langle w_{s_1, k}: k \leq j_n  \rangle$
is decreasing, we must have $n<i.$ It then follows that
$ h(n)< h(i)< n_i \leq j_n,$
and hence
\begin{center}
$(*)_2$ \quad\quad $r\Vdash$``$\name\tau_n=j_n >h(n)$''.
\end{center}
 By $(*)_1$ and $(*)_2$ we get a contradiction.
\end{proof}
It follows from Claim \ref{cla3} that the real $\name\eta[G_{\mathbb{Q}^1_{\bold n}}]$
dominates every ground model real and the lemma follows.
\end{proof}
We now prove Theorem \ref{t2}.
\begin{theorem}
\label{l-8}
Suppose $\bold n$ is a nice parameter. Then forcing with $\mathbb{Q}_{\bold n}^1$ adds a Hechler
real.
\end{theorem}
\begin{proof}
Let
$\mathbb D$ denote the Hechler forcing. Given $I\subseteq \mathbb D$ and
$f\in \omega^{\omega}$, let us say that \emph{$f$ satisfies $I$}, if there
exists a condition $p \in I$ such that $t_p \lhd f$ and $f_p(n) \leq f(n)$
for every $n<\omega$.
The next lemma shows that it suffices to verify that forcing with  $\mathbb{Q}_{\bold n}^1$
adds a real which satisfies $I$, for every maximal antichain $I\subseteq \mathbb D$ in the ground model.
\begin{lemma}
\label{le-1}
Suppose $\name\rho$ is a  $\mathbb{Q}_{\bold n}^1$-name for a real, such that for every maximal antichain $I\subseteq \mathbb D$ from the ground model,
$\Vdash_{\mathbb{Q}_{\bold n}^1}$``$\name\rho$
satisfies $I$''. Then
$\rho=\name\rho[G_{\mathbb{Q}_{\bold n}^1}]$ is a Hechler real.
\end{lemma}
\begin{proof}
We have to show that the set
\[
G=\{ (t, f) \in \mathbb{D}: t\lhd \rho \text{~and~} (\forall \ell g(t)  \leq n<\omega) f(n) \leq \rho(n)        \}
\]
is a $\mathbb{D}$-generic filter over $V$. Thus suppose that $I\subseteq \mathbb D$ is a maximal antichain in $V$.
By our assumption, there exists $(t, f) \in I$ such that $t \lhd f$ and $f(n) \leq \rho(n)$
for every $n<\omega$. It then follows that $(t, f) \in G \cap I$, and hence $G \cap I \neq \emptyset.$
The result follows.
\end{proof}
We now introduce a $\mathbb{Q}_{\bold n}^1$-name $\name\rho$ as requested by the above lemma.
 Let $\langle \underset{\sim}{\tau_n} : n<\omega \rangle$
be as in proof of Lemma \ref{l-7}. We define the $\mathbb{Q}_{\bold n}^1$-names $\langle \name{l}_i: i<\omega      \rangle$,
$\langle \name{k}_i: i<\omega      \rangle$
and $\name{\rho}$ as follows:
\begin{enumerate}
\item For every $i<\omega$, let $\name{l}_i$ is such that
\[
\Vdash_{\mathbb{Q}_{\bold n}^1}\text{``} \name{l}_i=\max\{l : \name{\eta}^1_{\bold n} \restriction (\name{\tau}_i+1) \in w_{\name{\eta}^1_{\bold n} \restriction \name{\tau}_i,i+l}\}\text{''},
\]

\item The name $\name{k}_i$ is defined by induction on
$i$ such that
\[
\Vdash_{\mathbb{Q}_{\bold n}^1}\text{``}\name{k}_i=\min\{k>i: (\forall j<i) k > \name{k}_j \text{~and~} \name{l}_k > 1\}\text{''},
\]

\item $\Vdash_{\mathbb{Q}_{\bold n}^1}\text{``}\name{\rho}=\langle \name{\tau}_n + \name{l}_{\name{k}_n} : n<\omega \rangle \in \omega^{\omega}$''.
\end{enumerate}
\begin{lemma}
\label{le-2}
Let $I=\{(t_n,f_n) : n<\omega\} \subseteq \mathbb D$ be a maximal
antichain and let $p\in \mathbb{Q}_{\bold n}^1$. Then there exists a condition
$q\in \mathbb{Q}_{\bold n}^1$ such that $p\leq q$ and $q\Vdash_{\mathbb{Q}_{\bold n}^1}$``$\name\rho$
satisfies $I$''.
\end{lemma}
\begin{proof}
Set $p_1=p.$
Let $h'\in \omega^{\omega}$ be a function such that for every $n<\omega$,
\begin{itemize}
\item $f_n \leq^* h'$,
\item $n<h'(n)<h'(n+1)$.
\end{itemize}
Let also $h\in \omega^{\omega}$ be defined as $h(n)=h'(n)+1$.

By the proof of Lemma \ref{l-7},  there are $p_2$ and $n_1^*$ such that:
\begin{itemize}
\item $p_1 \leq p_2$,
\item $2< n_1^* \leq \ell g(\tr(p_2))$,
\item $p_2 \Vdash$``$n_1^* \leq l \Rightarrow h(l) \leq \name\tau_l$''.
\end{itemize}
\begin{claim}
\label{cla4}
Assume $p \in \mathbb{Q}_{\bold n}^1$ and $\ell g(\tr(p))>2$. Then there exists an
increasing sequence $\langle n_i: i<\omega     \rangle$
of natural numbers satisfying the following conditions:
\begin{enumerate}
\item $n_0=\ell g(\tr(p))$,
\item If $l_1 \in [n_i,n_{i+1})$ and $t_1 \in T_p \cap \omega^{l_1}$,
then there are $l_2 \in [n_{i+1},n_{i+2})$ and $t_2 \in T_p \cap \omega^{l_2}$
such that:
\begin{enumerate}
\item $t_2$ extends $t_1$,
\item For every $l\in [l_1,l_2)$ we have
$t_2 \restriction (l+1) \notin w_{t_2 \restriction l,0}$,
\item $\beth_{i+1}(0)<\nor_{t_2}(\Suc_{T_p}(t_2))$.\footnote{$\beth_i(k)$ is defined by induction on $i$ by $\beth_0(k)=k$ and $\beth_{i+1}(k)=2^{\beth_i(k)}$.}
\end{enumerate}
\end{enumerate}
\end{claim}
\begin{proof}
By extending $p$ if necessary, assume that
 $\nor_{t}(\Suc_{T_p}(t))>2$ for every $\tr(p) \unlhd t \in T_p$.
Set $n_0=\ell g(\tr(p))$. Now
suppose
that $i<\omega$ and $n_{i+1}$ is defined. We define $n_{i+2}$.

Let $l_1 \in [n_i,n_{i+1})$ and $t_1 \in T_p \cap \omega^{l_1}$.  We define a sequence
$\langle s^{l_1, t_1}_{l} : l_1 \leq l<\omega \rangle$
 by induction on $l \geq l_1$
such that:
\begin{itemize}
\item $s^{l_1, t_1}_{l_1}=t_1$,
\item $s^{l_1, t_1}_{l+1} \in \Suc_{T_p}(s^{l_1, t_1}_l) \setminus w_{s^{l_1, t_1}_l,0}$.
\end{itemize}
We can easily define such a sequence, as for each $l$, $\Suc_{T_p}(s^{l_1, t_1}_l) \setminus w_{s^{l_1, t_1}_l,0} \neq \emptyset.$
 Let $\eta^{l_1, t_1}=\bigcup_{l_1 \leq l<\omega}s^{l_1, t_1}_{l}$.
then
\[
\lim_{n<\omega}(\nor_{\eta^{l_1, t_1} \restriction n}(\Suc_{T_p}(\eta^{l_1, t_1} \restriction n)))=\infty,
\]
and therefore there exists $n^{l_1, t_1}_{i+2} \geq n_{i+1}$ such that
$\nor_{\eta^{l_1, t_1} \restriction m}(\Suc_{T_p}(\eta^{l_1, t_1} \restriction m))>\beth_{i+1}(0)$
for every $n^{l_1, t_1}_{i+2} \leq m$.
Let
\[
n_{i+2} =\max\{ n^{l_1, t_1}_{i+2}: l_1 \in [n_i,n_{i+1}) \text{~and~} t_1 \in T_p \cap \omega^{l}    \}+1.
\]
Note that $n_{i+2}$ is well-defined as for each $l_1 \in [n_i,n_{i+1})$, the set $T_p \cap \omega^{l}$
is finite.
It is easy to see
that $n_{i+2}$ is as required.
\end{proof}
Let  $\langle n_i : i<\omega \rangle$ be a sequence as in Claim \ref{cla4} for the condition $p_2.$
Let $j_*$ and $j_{**}$ be such that:
\begin{itemize}
\item $j_*$ is maximal such that $p_2$ decides $\langle \name\tau_i: i < {j}_*        \rangle$,

\item $j_{**}$ is maximal such that $p_2$ decides $\langle \name k_i: i < {j_{**}}        \rangle$.
\end{itemize}
Let also the sequences $\langle m_i: i < j_*      \rangle$ and $\langle  k_i: i < j_{**}        \rangle$
be such that
\[
p_2 \Vdash \text{``} \langle \name\tau_i: i < {j_*}    \rangle= \langle m_i: i < j_*    \rangle    \text{''},
\]
and
\[
 p_2 \Vdash \text{``}  \langle \name k_i: i < {j_{**}}    \rangle=  \langle  k_i: i < j_{**}   \rangle  \text{''}.
\]

Let $s_1=\langle m_i+l_{k_i} : i<j_{**} \rangle$. Then
\[
p_2 \Vdash \text{``} \name\rho \restriction j_{**}=s_1 \text{''}.
\]
Define
$h_1 \in \omega^{\omega}$ by
$$h_1(i)=\beth_{h(i+1)+n_{i+2}}(0)+\max \{m_j : j<j_*\}+h(i+2).$$
Then $h_1 \geq h$ and  $(s_1, s_1 \cup h_1 \restriction [\ell g(s_1),\omega)) \in \mathbb D.$
As $I$ is a maximal antichain, there exists a condition  $(s_2, h_2) \in \mathbb D$
which extends $(s_1, s_1 \cup h_1 \restriction [\ell g(s_1),\omega))$
and a member of $I$.
We show that there is an extension of $p_2$ which forces ``$\name\rho$
satisfies $(s_2, h_2)$''.
\begin{claim}
\label{cla5}
There exists an increasing sequence $\langle p_{3, i}: i \in [\ell g(s_1), \ell g(s_2)]         \rangle$
of elements of $\mathbb{Q}^1_{\bold n}$ satisfying the following conditions:

\begin{enumerate}
\item For each $i \in [\ell g(s_1), \ell g(s_2))$, $p_{3, i}= (t^*_i, T_{p_2}^{[t^*_i]})$,
for some $t^*_i$.

\item $p_{3, \ell g(s_1)}=p_2$,
\item $p_{3,i}$ decides $\langle \name\tau_j: j < j_*+i-\ell g(s_1) \rangle$,
say
\[
p_{3, i}  \Vdash \text{``} \langle \name\tau_j: j < j_*+i-\ell g(s_1) \rangle =\langle m_j: j < j_*+i-\ell g(s_1) \rangle \text{''},
\]
\item $p_{3,i}$ doesn't force a value for $\name\tau_{j_*+i-\ell g(s_1))}$.
\end{enumerate}
\end{claim}
\begin{proof}
We define a sequence $\langle t^*_i: i \in [\ell g(s_1), \ell g(s_2)]         \rangle$,
by induction on $i$ as follows. Let $t^*_{\ell g(s_1)}=\tr(p_2).$
Now suppose that $t^*_i$ is defined and $\ell g(t^*_i) \in [n_{i-\ell g(s_1)}, n_{i+1-\ell g(s_1)}]$.
By Claim \ref{cla4}, applied to $\ell g(t^*_i)$ and $t^*_i$, there exist
$l^* \in [n_{i+1-\ell g(s_1)}, n_{i+2-\ell g(s_1)}]$
and $t^* \in T_{p_2} \cap \omega^{l^*}$ such that:
\begin{itemize}
\item $t^*_i \unlhd t^*$ ,
\item For every $l \in [\ell g(t^*_i), l^*), t^* \restriction (l+1) \notin w_{t^* \upharpoonright l, 0}$,
\item $\beth_{i-\ell g(s_1)+1}(0)<\nor_{t^*}(\Suc_{T_p}(t^*))$.
\end{itemize}
Let $u=\Suc_{T_{\bold n}}(t^*) \setminus \Suc_{T_{p_2}}(t^*)$ and $j=j_*+i - \ell g(s_1)$. We show that
$w_{t^*,j} \setminus w_{t^*,j+1} \nsubseteq u$. Suppose not. It then follows that
$$|w_{t^*,j+1}| \leq \frac{|w_{t^*,j}|}{2} \leq |w_{t^*,j} \setminus w_{t^*,j+1}| \leq |u|,$$
and
therefore,
by the
construction of condition $p_2$ and the choice of the $w$'s,
$$\beth_{i-\ell g(s_1)+1}(0) < \nor_{t^*}(\Suc_{T_{p_2}}(t^*))=\nor_{t^*}(\Suc_{T_{\bold n}}(t^*) \setminus u) \leq i-\ell g(s_1)+2,$$
which is not possible.
Therefore $w_{t^*,j} \setminus w_{t^*,j+1} \nsubseteq u$, and hence we can find some
$t^*_{i+1} \in \Suc_{T_{p_2}}(t^*) \cap (w_{t^*, j} \setminus w_{t*,j+1}) \neq \emptyset.$
This defines $t^*_{i+1}.$

It is now easy to see that $p_{3, i}= (t^*_i, T_{p_2}^{[t^*_i]})$'s are as required.
\end{proof}

Note that each $p_{3, i}$ from the above claim also decides  $(\name l_j: j<j_*+i-\ell g(s_1) \rangle,$ say
\[
p_{3, i}  \Vdash \text{``} \langle \name l_j: j < j_*+i-\ell g(s_1) \rangle =\langle l^*_j: j < j_*+i-\ell g(s_1) \rangle \text{''}.
\]
Note also that for every $i\in (\ell g(s_1), \ell g(s_2))$,
$$m_{i-1} \leq \ell g(\tr(p_{3,i}))=\ell g(t^*_i) \leq n_{i+1} \leq \beth_{h(i)+n_{i+1}}(0) < h_1(i-1) \leq  h_2(i-1).$$
Now choose $p_4$ such that:
\begin{itemize}
\item $p_{3, \ell g(s_2)} \leq p_4$,
\item $\max(\rng(s_2))<\ell g(\tr(p_4))$,
\item $p_4$ does not force a value for $\name\tau_{j_*+\ell g(s_2)-\ell g(s_1)}$\footnote{This can be done easily, for example, by extending $\tr(p_{3,\ell g(s_2)})$
at each stage to a sequence outside of the appropriate $w_{s,0}$}.
\end{itemize}
\begin{claim}
\label{cla6}
There exists an increasing sequence $\langle p_{5, i}: i \in [\ell g(s_1), \ell g(s_2)]       \rangle$
of elements of $\mathbb{Q}^1_{\bold n}$ satisfying the following conditions:
\begin{enumerate}
\item $p_{5,\ell g(s_1)}=p_4$,

\item  $p_{5,i}$ forces a value for $\name\tau_{j_*+(\ell g(s_2)-\ell g(s_1))+(j-\ell g(s_1))}$
iff $\ell g(s_1) \leq j<i$. In this case, let $m_{j_*+(\ell g(s_2)-\ell g(s_1))+(j-\ell g(s_1))}$ be such that
\[
p_{5, i} \Vdash \text{``} \name\tau_{j_*+(\ell g(s_2)-\ell g(s_1))+(j-\ell g(s_1))}=   m_{j_*+(\ell g(s_2)-\ell g(s_1))+(j-\ell g(s_1))}      \text{''},
\]
\item For $\ell g(s_1) \leq j<i$,
\[
p_{5, i} \Vdash \text{``}\name{k}_{j}=j_*+(\ell g(s_2)-\ell g(s_1))+(j-\ell g(s_1))\text{'',}
\]
\item For $\ell g(s_1) \leq j<i$,
\[
p_{5, i} \Vdash \text{``}\name{l}_{j_*+(\ell g(s_2)-\ell g(s_1))+(j-\ell g(s_1))}=s_2(j)-m_{j}\text{''.}
\]
\end{enumerate}
\end{claim}
\begin{proof}
We define $p_{5, i}$'s, by induction on $i$. Set $p_{5,\ell g(s_1)}=p_4$.
Now suppose that $i, i+1 \in [\ell g(s_1), \ell g(s_2)]$
and $p_{5, i}$ is defined.

Set $s_i^*=\tr(p_{5,i})$ and $i^*=s_2(i)- m_{i}$. By the choice of the sequence $\langle n_i: i<\omega   \rangle$,
we have $1 <i^*.$
Also, for notational simplicity set
$$j_i=j_*+(\ell g(s_2)-\ell g(s_1))+(i-\ell g(s_1)).$$
 By an argument similar to the proof of Claim \ref{cla4},
we can find $s^*$ such that:
\begin{itemize}
\item $s_i^* \unlhd s^*,$

\item For
every $l\in (\ell g(s_i^*), \ell g(s^*))$, $s^* \restriction (l+1) \notin w_{s^* \restriction l,0}$,

\item  $j_i+i^*+2< \nor_{s^*}(\Suc_{T_{p_{5, i}}}(s^*))$.
\end{itemize}
Choose
$$s^*_{i+1} \in \Suc_{T_{p_{5,i}}}(s^*) \cap \big(w_{s^*, j_i+i^*} \setminus w_{s^*, j_i+i^*+1}\big)$$
and define $p_{5,i+1}:=(s^*_{i+1}, T_{p_{5, i}}^{[s^*_{i+1}]})$.
Then $p_{5, i+1} \geq p_{5, i}.$
 We show that it satisfies items (2)-(4) of the claim.

{\bf $p_{5, i+1}$ satisfies clause (2).} To show this, we need to have:
\begin{itemize}
\item $\tr(p_{5,i+1})=s^*_{i+1} \in w_{s^*,j_i}$,
\item Every initial segment $t$ of $s^*_{i+1}$
avoids $w_{t \restriction (\ell g(t)-1),j_i}$.
\end{itemize}
These are clearly true by the way we defined $s^*_{i+1}$.

{\bf $p_{5, i+1}$ satisfies clause (3).} By the choice of the conditions $p_{3, j}$,
$l_j \leq 1$ for every $j<j_*+(\ell g(s_2)-\ell g(s_1))$. Since $p_{5, i+1}$
extends $p_{5, i}$, it follows from  the induction hypothesis that $p_{5, i+1}\Vdash$``$\name{k}_{j}=j_j$'' for all $j<i$.
Thus to guarantee clause (3),
we only have to show that $p_{5, i+1}\Vdash$``$\name{k}_{i}=j_i$''.
By the choice of the sequence $\langle n_i: i<\omega \rangle$ and the conditions
$p_{3, i}$ we have $i^* \geq 2$. Furthermore $p_{5, i+1} \Vdash\text{``}\name\tau_{j_i}=m_{j_i}  $'' and
\begin{center}
$p_{5, i+1} \Vdash\text{``}  \name\eta^1_{\bold n} \upharpoonright  (m_{j_i} +1) \in w_{\name\eta^1_{\bold n}  \restriction m_{j_i} , j_i+i^*}  \text{''}$.
 \end{center}
 In particular, it follows that
$p_{5, i+1} \Vdash\text{``} \name{l}_{j_i} \geq i^* > 1$''. Thus by its definition
\[
p_{5, i+1} \Vdash\text{``}  \name{k}_{i} =\min \{k >i: (\forall j<i) k > \name{k}_j \text{~and~} \name l_k > 1 \} = j_i                            \text{''},
\]
which gives the result.

{\bf $p_{5, i+1}$ satisfies clause (4).} This is clear, as $s^*_{i+1} \in w_{s^*, j_i+i^*} \setminus w_{s^*, j_i+i^*+1}$.
\end{proof}
Let $p_5=p_{5,\ell g(s_2)}$. It's easy to see  that
$p_5 \Vdash$``$s_2 \unlhd \name\rho$''. Indeed, as $p_2 \leq p_5,$ we have $p_5 \Vdash$``$\name\rho \upharpoonright \ell g(s_1)=s_1=s_2 \upharpoonright \ell g(s_1)$''. On the other hand, for every
$\ell g(s_1) \leq j < \ell g(s_2)$ we have, using Claim \ref{cla6},
\[
p_5 \Vdash \text{``} s_2(j) =\name{l}_{\name{k}_j}+m_j= l_{j_*+(\ell g(s_2)-\ell g(s_1))+(j-\ell g(s_1))}+m_j=\name\rho(j)\text{''}.
\]
\begin{claim}
\label{cla8}
There exists a condition $p_6$ such that $p_5 \leq p_6$ and
 $p_6 \Vdash$``$h_2(l) \leq \name\rho(l)$''
for every $\ell g(s_2) \leq l$.
\end{claim}
\begin{proof}
 We already know, by the choice of the condition $p_2$,
that there exists $m_*$ such that for every $l > m_*$,
$p_2$ an hence $p_5$ forces ``$h_2(l) \leq \name\rho(l)$''.

 By the proof of Lemma \ref{l-7}, there is a condition $p_5 \leq p'$ such that
$\tr(p')=\tr(p_5)$ and $p' \Vdash$``$\ell g(\tr(p_5)) \leq n \Rightarrow h_2(n) \leq \name\tau_n$''.
Therefore, we may assume  that $p_5$ decides $\name\tau_n$
for every $n\in [\ell g(s_2),m_*]$.

On the other hand, $p_5$ does not
decide $\name{l}_{\name{k}_n}$,  for  $n\in [\ell g(s_2),m_*]$, as the trunk of $p_5$
is the first place where $\name{l}_{\name{k}_{\ell g(s_2)-1}}$ is
decided.
Thus we can repeat the argument that lead us from $p_2$
to $p_5$ in order to find a condition $p_6$ such that:
\begin{itemize}
\item $p'\leq p_6,$
\item $p_6 \Vdash$``$h_2(n) \leq \name\rho(n)$
for every $n\in [\ell g(s_2),m_*]$.
\end{itemize}
Then $p_6$ is as required.
\end{proof}
Finally let $q=p_6.$ Then $q$ extends $p$ and it forces ``$\name\rho$ satisfies $I$''.
This completes the proof of Lemma \ref{le-2}.
\end{proof}
By Lemma \ref{le-2}, $\name\rho[G_{\mathbb{Q}^1_{\bold n}}]$ satisfies $I$, for every maximal antichain
$I \subseteq \mathbb D$ in the ground model, hence by Lemma \ref{le-1}, it is a Hechler real.
Theorem \ref{l-8} follows.
\end{proof}

\begin{corollary} ($\ZF+\DC_{\omega_1}$)
\label{c-22}
 Let $\bold n$ be
 a nice parameter.
Suppose every set of reals is $I_{\mathbb Q^1_{\bold n},\aleph_0}$-measurable.
Then there is an inner model of $\ZFC$ with a measurable cardinal.
\end{corollary}
\begin{proof}
By Theorem \ref{l-8}, forcing with $\mathbb Q^1_{\bold n}$
adds a Hechler real. Now the result follows from Corollary \ref{c-2}.
\end{proof}
We close the paper by the following question, which asks
whether a measurable cardinal is an optimal lower bound on the consistency strength of the results obtained above.
\begin{question}
What is the consistency strength of ``$\ZF+\DC_{\omega_1}+$every set of reals is $I_{\mathbb Q,\aleph_0}$-measurable'',
where $\mathbb Q$ is a Suslin ccc forcing
notion which adds a Hechler real.
\end{question}

\end{document}